\documentclass[11pt,leqno]{amsart}
\PassOptionsToPackage{dvipsnames}{xcolor}

\usepackage[letterpaper,margin=1.17in,footskip=0.27in]{geometry}
\usepackage[T1]{fontenc}
\IfFileExists{scholax.sty}{%
  \usepackage[p,scosf,looser]{scholax}
  \usepackage[scaled=1.08,ncf,vvarbb]{newtxmath}

  \usepackage[cal=stixtwoplain]{mathalfa}
}{%
  \usepackage{lmodern}
  \usepackage{mathpazo}
  \usepackage{amssymb}
}

\DeclareSymbolFont{opnames}{T1}{\rmdefault}{m}{n}
\SetSymbolFont{opnames}{bold}{T1}{\rmdefault}{b}{n}
\makeatletter\renewcommand{\operator@font}{\mathgroup\symopnames}\makeatother
\usepackage{amsmath}
\usepackage{amsthm}
\usepackage{stmaryrd}
\usepackage{mathrsfs}
\usepackage{mathtools}
\usepackage{aliascnt}
\usepackage{microtype}
\usepackage{graphicx}
\usepackage{verbatim}
\usepackage{xcolor}

\usepackage{tikz-cd}
\usepackage{enumitem}
\usepackage[colorlinks,pagebackref]{hyperref}

\usepackage[capitalise,nameinlink,noabbrev]{cleveref}
\usepackage{bookmark}
\usepackage{todonotes}
\crefformat{equation}{#2(#1)#3}
\Crefformat{equation}{#2(#1)#3}
\crefrangeformat{equation}{#3(#1)#4--#5(#2)#6}
\crefmultiformat{equation}{#2(#1)#3}{ and~#2(#1)#3}{, #2(#1)#3}{, and~#2(#1)#3}

\hypersetup{citecolor=darkforestgreen,linkcolor=Mahogany,urlcolor=darkforestgreen,
  pdftitle={Uniform Positivity of the F-Signature}}

\definecolor{darkmaroon}{RGB}{110,0,20}
\definecolor{darkforestgreen}{RGB}{10,80,30}
\definecolor{Burgundy}{RGB}{144,0,32}
\definecolor{GreenH}{rgb}{0.0, 0.5, 0.0}
\definecolor{calpolypomonagreen}{rgb}{0.12, 0.3, 0.17}
\definecolor{utahcrimson}{rgb}{0.83, 0.0, 0.25}
\definecolor{umber}{rgb}{0.39, 0.32, 0.28}
\definecolor{regalia}{rgb}{0.32, 0.18, 0.5}

\allowdisplaybreaks
\numberwithin{equation}{section}
\newtheorem{theorem}{Theorem}[section]
\newtheorem{lemma}[theorem]{Lemma}
\newtheorem{proposition}[theorem]{Proposition}
\newtheorem{conjecture}[theorem]{Conjecture}
\newtheorem{corollary}[theorem]{Corollary}

\theoremstyle{definition}
\newtheorem{definition}[theorem]{Definition}

\newtheorem{notation}[theorem]{Notation}

\theoremstyle{remark}
\newtheorem{remark}[theorem]{Remark}

\crefname{notation}{Notation}{Notations}
\crefname{assumption}{Assumption}{Assumptions}
\crefname{question}{Question}{Questions}
\crefname{problem}{Problem}{Problems}
\crefname{conjecture}{Conjecture}{Conjectures}
\newcommand{\introthmname}{}
\newtheorem{introthminn}{\introthmname}

  \newcommand{\introdefname}{}
  \newtheorem*{introdefinn}{\introdefname}

\renewcommand{\O}{\mathcal{O}}

\newcommand{\C}{\mathbb{C}}

\newcommand{\N}{\mathbb{N}}
\newcommand{\Q}{\mathbb{Q}}

\newcommand{\Z}{\mathbb{Z}}

\newcommand{\vare}{\varepsilon}
\newcommand{\Spec}{\operatorname{Spec}}

\newcommand{\Id}{\operatorname{id}}

\newcommand{\G}{\mathbb{G}}

\newcommand{\lre}{\longrightarrow}

\newcommand{\Div}{\operatorname{Div}}

\newcommand{\vol}{\operatorname{vol}}

\newcommand{\ord}{\operatorname{ord}}
\providecommand{\frm}{}
\renewcommand{\frm}{\mathfrak{m}}

\newcommand{\llb}{\llbracket}
\newcommand{\rrb}{\rrbracket}
\newcommand{\Supp}{\operatorname{Supp}}

\newcommand{\stacktag}[1]{\href{https://stacks.math.columbia.edu/tag/#1}{#1}}
\DeclareMathOperator{\Hom}{Hom}
\DeclareMathOperator{\Frac}{Frac}

\title{Uniform Positivity of the $F$-Signature and Pure Descent of klt Singularities}
\author{Rahul Ajit}
\address{Department of Mathematics, The University of Kansas, 405 Snow Hall, Lawrence, KS 66045}
\email{rahulajit@ku.edu}
\author{Daniel Apsley}
\address{Department of Mathematics, University of Utah, Salt Lake City, UT, 84012, USA.}
\email{apsley@math.utah.edu}
\date{\today}

\begin{document}

\begin{abstract}
In this paper, we affirmatively answer the question of Carvajal-Rojas, Schwede, and Tucker on the uniform positivity of the $F$-signature, and prove a converse. Namely, given a klt type singularity, we establish uniform lower bounds on the $F$-signatures of its reductions mod $p$ for $p$ sufficiently large. Conversely, if the $F$-signatures of the reductions of a normal singularity are uniformly bounded below by a positive constant on large primes, we prove that the singularity is of klt type. As an application, we give a new reduction-mod-$p$ proof of Zhuang's theorem that klt type singularities descend along pure morphisms.
\end{abstract}

\maketitle

\section{Introduction}
The past few decades have seen spectacular progress towards understanding the relationship between singularities appearing in the minimal model program and analogous notions in positive characteristic, defined using the Frobenius map. \textit{Kawamata log terminal (klt) singularities}, defined via discrepancies of exceptional divisors on a log resolution, are a central class of singularities in the minimal model program (see \cite[Definition 2.34]{kollarmori}). If $x \in X$ is a closed point of a normal complex variety and $\Delta$ is an effective $\Q$-divisor on $X$ so that $K_X + \Delta$ is $\Q$-Cartier near $x$, then it follows from \cite{HW02,takagi04} that $x \in (X, \Delta)$ is klt if and only if for large primes $p$, the reductions $x_p \in (X_p, \Delta_p)$ are strongly $F$-regular. Unlike the definition of klt singularities, the definition of strong $F$-regularity does not require the singularity to be $\Q$-Gorenstein. If $x \in X$ is \textit{of klt type}, meaning that there exists a boundary $\Delta$ such that $x \in (X, \Delta)$ is a klt singularity, then the reductions mod $p$ of the singularity $x \in X$ itself are strongly $F$-regular for large $p$. The converse is a long-standing open problem. 

\begin{conjecture}[{\cite[Problem 5.3.1]{HW02}}] \label{conj: HWGFRtype}
Let $x \in X$ be a closed point of a normal complex variety. Then $x \in X$ is of klt type if and only if for large $p$, the reductions $x_p \in X_p$ mod $p$ are strongly $F$-regular. 
\end{conjecture}

We do not resolve this conjecture in this article, of course, but we give a different characterization of klt type singularities by reduction mod $p$, via the \textit{$F$-signature}. For an $F$-finite local ring $R$ of dimension $n$ and characteristic $p > 0$ with perfect residue field, the $F$-signature, introduced in \cite{ hl02}, is given by the limit  $$s(R) = \lim_{e \to \infty} \frac{\operatorname{frk}_R(F^e_* R)}{p^{en}},$$ which exists in this generality by \cite{tucker12}. The $F$-signature characterizes strong $F$-regularity by \cite{AL03}; $R$ is strongly $F$-regular if and only if $s(R) > 0$.

In this article, we prove that in the klt type setting, the $F$-signatures of reductions are uniformly bounded below by a positive constant, thus affirmatively answering \cite[Question 5.9]{crst18} (see also \cite[Conjecture 1.1]{lp26} and \cite[Conjecture 1.2]{TY25}). Furthermore, we prove a converse: uniform positivity of the $F$-signature characterizes klt type singularities.

\subsection*{Main Theorem} To state the theorem, we spread out $x \in X$ to a \textit{model}: a finitely generated $\Z$-algebra $\Lambda \subseteq \C$, a finite-type $\Lambda$-scheme $X_\Lambda$, and a $\Lambda$-valued point $x_\Lambda \in X_\Lambda(\Lambda)$, so that the base change $x_\C \in X_\C$ is isomorphic to $x \in X$.

\begin{theorem} \label{thm: main}
 Let $\Lambda \subset \C$ be a finite type $\Z$-algebra and let $X_\Lambda \to \Spec \Lambda$ be a finite-type $\Lambda$-scheme. Assume $x_\Lambda \in X_\Lambda(\Lambda)$ is a $\Lambda$-valued point, and set $x \in X$ to be the base change $x_\C \in X_\C$. Suppose $X$ is normal near $x$. 
\begin{enumerate}[ref=\thetheorem(\arabic*)]
\item\label[theorem]{thm: main1}(see \cref{thm: forwardirection}) If $x \in X$ is a klt type singularity, then there exist a real number $\vare > 0$ and a nonempty open subset $U \subseteq \Spec \Lambda$ so that $s(\O_{X_b, x_b}) \geq \vare$ for all closed points $b$ in $U$.
\item\label[theorem]{thm: main2}(see \cref{thm: conversedirection}) Conversely, if $x \in X$ is not of klt type, and $\vare > 0$ is a real number, then there is a nonempty open subset $U \subseteq \Spec \Lambda$ with $s(\O_{X_b, x_b}) < \vare$ for all closed points $b$ in $U$.
\end{enumerate}
\end{theorem}

In particular, if $\Lambda = \Z$, then \cref{thm: main} implies that $x \in X$ is a klt type singularity if and only if $\liminf_{p \to \infty} s(\O_{X_p, x_p}) > 0$ if and only if $\limsup_{p \to \infty} s(\O_{X_p, x_p}) > 0$. Moreover, with this theorem, settling \cref{conj: HWGFRtype} amounts to determining whether or not strongly $F$-regular type singularities admit uniform lower bounds for the $F$-signature.

To prove the forward direction of \cref{thm: main}, we spread out a boundary associated to a Noether normalization, similar to \cite[Proposition 6.4]{bgo17}. After spreading out the necessary data to a model and reducing mod $p$, this direction rests on the following lower bound on the $F$-signature of a certain pair. The $F$-signature was extended to pairs in \cite{bst12}, where for an effective $\Q$-divisor $\Delta$ on $\Spec R$, we have $s(R, \Delta) \leq s(R)$.

\begin{theorem}[{cf. \cref{thm:charpmain}}]
    Suppose $(R, \frak{m})$ is a local normal domain of dimension $n$, essentially of finite type over a perfect field $k = R/\frm$ of characteristic $p > 0$. Fix an effective $\Q$-divisor $\Delta$ on $\Spec R$ with denominators prime to $p$. In this setting, there exists a nonzero $d \in R$, depending on a Noether normalization, satisfying the following. 
    
    If $\eta \in \Q\cap (0,1)$ has denominator prime to $p$ and $(R, \Delta + \eta \Div(d))$ is strongly $F$-regular, then $$s(R, \Delta) \geq \frac{\eta^n}{2^n n!} > 0.$$
\end{theorem}

The above is obtained using multi-variable \textit{Hasse--Schmidt derivatives} to produce enough splittings. The converse \cref{thm: main2} is obtained by relating splittings to divisors with positive log discrepancies in a uniform way.

Prior to this work, \cref{thm: main1} was known for pure subrings of regular local rings, and it was also known that this direction reduces to the case of Gorenstein singularities; see \cite[Corollary 1.4]{TY25}. This direction was also proven for non-weakly exceptional threefold singularities in \cite[Theorem 1.5]{lp26}.

\subsection*{Applications} 
As was noticed in \cite{crst18}, the uniform positivity of the $F$-signature has a number of notable consequences. 

\subsubsection*{Global Analogue} By fixing an ample divisor $A$ on a normal projective variety $Y$, we obtain a reduction-mod-$p$ criterion for being of Fano type by considering the affine cone $C(Y, A) = \Spec \bigoplus_{m \geq 0} H^0(Y, \O_Y(mA))$ over $Y$. 

\begin{corollary} (see \cref{cor: section_global_version}) \label{cor: global_version}
Let $Y$ be a normal projective variety and $A$ be an ample Cartier divisor on $Y$. Then, $Y$ is of Fano type if and only if all models of the cone point $o \in C(Y, A)$ satisfy the conclusion of \cref{thm: main1}.
\end{corollary}

Using \cite[Lemma 3.7]{lp26}, we note that the conclusion of \cref{thm: main1} applied to $o \in C(Y, A)$ is equivalent to the condition that the Frobenius $\alpha$-invariants (cf. \cite{pande2025frobeniusversiontiansalphainvariant})  of the reductions $(Y_p, A_p)$ are uniformly bounded below by a positive constant. Thus, we obtain a positive answer to \cite[Conjecture 1.2]{lp26} as well. 

\subsubsection*{Pure Descent} We are also able to adapt the proof of \cref{thm: main} to give a reduction-mod-$p$ proof of Zhuang's theorem \cite[Theorem 1.1]{zhuang24} on the pure descent of klt type singularities. (See also \cite[Theorem B]{Schoutens05}, \cite[Theorem 1]{BGLM24}, \cite[Theorem 1.3]{TY25}.) Recall that a ring homomorphism $A \to B$ is \textit{pure} if $M \to M \otimes_A B$ is injective for every $A$-module $M$. Examples of pure maps include affine GIT quotients and direct summands.

\begin{theorem}[{\cite[Theorem 1.1]{zhuang24}}] \label{thm: puredescent}
Let $f: \Spec B \to \Spec A$ be a pure morphism of normal affine complex varieties. If $\Spec B$ is of klt type, then so is $\Spec A$.
\end{theorem}

An obstacle to reducing this statement mod $p$ is that, unless $A \to B$ is finite, its reductions mod $p$ may not be pure; see \cite{HJSP23}. We get around this by constructing a finitely generated $A$-submodule of $B$ which, after reduction mod $p$, contains the values on $A$ of a generating set of $p^{-e}$-linear maps on $B$, so that purity need only be used on this submodule. In particular, our proof does not use the minimal model program.

\subsubsection*{Fundamental Groups} Lastly, combining \cref{thm: main} with \cite[Corollary 3.5]{crst18}, we obtain effective bounds, comparable to \cite[Theorem B]{JS22} and \cite[Corollary 1.4]{xz21}, on the orders of the \'etale fundamental groups of klt type singularities, as anticipated in \cite[Section 5]{crst18}. This gives an effective version of the finiteness results of \cite{xu14, bgo17}.

\begin{corollary} \label{cor: localfundgroups} If $x \in X$ is a klt type singularity over $\C$, then its local \'etale fundamental group is finite, of order at most $1/\varepsilon$, where $\varepsilon > 0$ is as in \cref{thm: main1}.
\end{corollary}

\subsection*{Acknowledgments}
The authors are very grateful to Harold Blum for carefully reading and commenting on earlier drafts. We would also like to thank Karl Schwede for informing us of the independent and overlapping work of \cite{LuoPolstra26}, and for his comments on earlier drafts. These comments have significantly improved our presentation. We would also like to thank Tommaso de Fernex  for suggesting to study Hasse--Schmidt and Frobenius together. Lastly, we thank Thomas Polstra for discussing our overlapping works in detail during the final stage of the writing process. 

\subsection*{AI Statement}   We collaborated extensively with Claude Opus 5.5 (Anthropic) at the start of this project. It began with \cite[Problem 46]{xu2026openproblemskstabilityfano}, which asks whether $\widehat{\vol}(x, X, \Delta) = \lim_{p \to \infty} \operatorname{Fvol}(x, X_p, \Delta_p)$. Prompted whether the Hasse--Schmidt methods of \cite{rahulmatthewbertini} could address this problem, it did not resolve it, but it proved a version of \cref{cor: splittingidealinclusions}, namely that splitting ideals are invariant under modified Hasse--Schmidt operators. Further human work with Claude led to \cref{thm:charpmain}, and hence we deduced the forward direction of \cref{thm: main}, and to \cref{prop: discrepancybound}, from which we obtained the converse. Later interaction with Claude produced a version of \cref{lem: finitemodule}, from which we were able to complete the proof of \cref{thm: puredescent}. We made several modifications and simplifications to Claude's inputs as well, based on our human taste. Claude was also used to identify gaps in earlier drafts, locate references, and edit the exposition. The authors have checked every argument and take full responsibility for the contents of this paper. All contents of this paper are human-written. 

\subsection*{Overlap with Independent Work} Luo and Polstra \cite{LuoPolstra26} independently established \cref{thm: main1} (see \cite[Main Theorem 2]{LuoPolstra26}). Their proof is related to ours; compare \cite[Main Theorem 1]{LuoPolstra26} to our \cref{thm:charpmain}. Both approaches rely on spreading out separable Noether normalizations and using formal Taylor expansions and clearing denominators via powers of a distinguished element. (Compare \cite[Lemma 3.1 and Corollary 3.2]{LuoPolstra26} to \cref{lem: Psi} and \cref{cor: splittingidealinclusions} below.) Their distinguished element is obtained by a Jacobian determinant, and their lower bound differs from ours in that it involves the Hilbert--Samuel multiplicity and an Izumi--Rees constant. 

Not in \cite{LuoPolstra26} is the converse statement \cref{thm: main2}, nor the corollaries: reduction-mod-$p$ criteria for Fano type (see \cref{cor: global_version}) and a reduction-mod-$p$ proof of the pure descent of klt type singularities (see \cref{thm: puredescent}). However, \cite{LuoPolstra26} addresses a number of interesting problems and applications which are not considered in this article. They give precise quantitative relationships between local $\alpha_F$-invariants, multiplicities, and valuations (see \cite[Main Theorems 1 and 4]{LuoPolstra26}). They also deduce lower bounds on the limit $F$-signature involving the normalized volume (see \cite[Main Theorem 3]{LuoPolstra26}). Lastly, they give an effective solution to the uniform symbolic topology problem for strongly $F$-regular local rings in positive characteristic, and for klt type local rings in characteristic zero (see \cite[Main Theorem 4]{LuoPolstra26}).

\section{Taylor Expansions and Uniform Lower Bounds}

In this section, we prove the characteristic $p > 0$ part of \cref{thm: main1}, given below as \cref{thm:charpmain}. The main idea will be to introduce Hasse--Schmidt operators, obtained via Taylor coefficients (see also \cite[16.11.1--16.11.2]{EGAIV}). Namely, we will prove that the splitting ideals are invariant under these operators in \cref{cor: splittingidealinclusions}. Using this invariance, we will bound the colength of these splitting ideals in \cref{prop: colengthbound} to prove the desired result. 

\subsection*{Hasse--Schmidt Operators} If $k$ is a field and $K = k(t_1, \dots, t_n)$, there is a $k$-algebra homomorphism $$\Theta_K: K \lre K \llb T_1, \dots, T_n \rrb$$ given by $f(t_i) \mapsto f(t_i + T_i)$, which maps $k[t_1, \dots, t_n]$ into $k[t_1, \dots, t_n][T_1, \dots, T_n]$. For $\alpha \in \N^n$, we denote the coefficient of $T^\alpha:= T^{\alpha_1}_1 \cdots T^{\alpha_n}_n$ in $\Theta_K(f)$ by $D^{[\alpha]}(f)$. Using the binomial theorem, these coefficients can explicitly be described for monomials as \begin{align} \label{eq:HSmonomial} D^{[\alpha]}(t^{\beta}) = \binom{\beta_1}{\alpha_1} \cdots \binom{\beta_n}{\alpha_n} t^{\beta - \alpha}. \end{align} Above, we take the convention that the binomial coefficients are zero when $\beta - \alpha \notin \N^n$, and that they are evaluated in $\Z$ before being mapped to $k$. In particular, if $n = 1$, we have $D^{[p]}(t^p) = 1$ for every prime $p$, regardless of the characteristic of $k$. Since $\Theta_K$ is a ring homomorphism, these coefficients satisfy the following Hasse--Schmidt Leibniz rule: \begin{align} \label{HSLeibniz} D^{[\alpha]}(fg) = \sum_{\beta + \gamma = \alpha} D^{[\beta]}(f) D^{[\gamma]}(g). \end{align} We call $D^{[\alpha]}: K \to K$ the \textit{Hasse--Schmidt operator of degree $\alpha$} on $K$. Using a multi-variable analogue of \cite[Theorem 1.5]{traves00} below, we can extend these operators to finite separable field extensions.

\begin{proposition} \label{prop:HSTaylorExpansion}
If $L$ is a finite separable extension of $K = k(t_1, \dots, t_n)$, there is a unique extension of $\Theta_K$ to a ring homomorphism $\Theta_L: L \to L \llb T_1, \dots, T_n \rrb$ satisfying $\Theta_L(f) \equiv f \bmod (T_1, \dots, T_n)$ for all $f \in L$.
\end{proposition}

\begin{proof}
Since $K \subset L$ is finite and separable, it is \'etale, so it is formally \'etale by \cite[Section~10.150, \href{https://stacks.math.columbia.edu/tag/00UP}{Tag~00UP}]{stacks-project}. Thus, given a suitable extension $$\Theta^m_L: L \lre L [T_1, \dots, T_n ]/(T_1, \dots, T_n)^m$$ of the reduction $\Theta^m_K$ of $\Theta_K$ modulo $(T_1, \dots, T_n)^m$, we can use formal \'etaleness to extend it to a unique $\Theta^{m + 1}_L$:
\[
\begin{tikzcd}
K \arrow[r, "\Theta^{m+1}_K"] \arrow[d] & L [T_1, \dots, T_n ]/(T_1, \dots, T_n)^{m+1} \arrow[d] \\
L \arrow[r, "\Theta^m_L"'] \arrow[ur, dashed, "\Theta^{m+1}_L" description] & L [T_1, \dots, T_n ]/(T_1, \dots, T_n)^m.
\end{tikzcd}
\]
Taking $\Theta_L^1= \Id_L$, these define the required extension $$\Theta_L: L \to \varprojlim_m L [T_1, \dots, T_n ]/(T_1, \dots, T_n)^m = L \llb T_1, \dots, T_n \rrb.$$
\end{proof}

\begin{definition} \label{def: HasseSchmidtOperators}
Suppose $L$ is a finite separable extension of $K$, and $\Theta = \Theta_L: L \to L\llb T_1, \dots, T_n \rrb$ is as above. The coefficients of $$\Theta(f) = \sum_{\alpha \in \N^n} D^{[\alpha]} (f) T^\alpha,$$ viewed as $k$-linear maps $D^{[\alpha]}: L \to L$, are called the \textit{Hasse--Schmidt operators of degree $\alpha$} on $L$.
\end{definition}

\begin{remark}
Since $\Theta$ is a ring homomorphism, these operators $D^{[\alpha]}: L \to L$ still satisfy the Hasse--Schmidt Leibniz rule \cref{HSLeibniz}.
\end{remark}

\begin{lemma}
   \label{rem: froblinearity}
   Suppose the characteristic of $k$ is $p > 0$, and let $e$ be a positive integer. If $\alpha \in \N^n$ satisfies $|\alpha| := \sum_j \alpha_j < q := p^e$, then $F^e_*D^{[\alpha]}: F^e_* L \to F^e_*L$ is $L$-linear.
\end{lemma}
\begin{proof}
 Using the fact that $\Theta$ is a ring homomorphism, as well as the congruence part of \cref{prop:HSTaylorExpansion}, we have that $$\Theta(f^q) = \Theta(f)^q \in f^q + (T_1^q, \dots, T_n^q)$$ in characteristic $p$. Thus, $D^{[\beta]}(f^q) = 0$ for $0 < |\beta| < q$. Applying the Hasse--Schmidt Leibniz rule \cref{HSLeibniz}, we then obtain $$D^{[\alpha]}(f^q g) = f^q D^{[\alpha]}(g)\text{ for all } f, g \in L,$$ as required.  
\end{proof}

To specialize to the setting of \cref{thm:charpmain} below, we fix the following notation. 

\begin{notation} \label{not:setup}
Let $(R, \frm)$ be a normal local domain of dimension $n$, essentially of finite type over a perfect field $k = R/\frm$ of characteristic $p > 0$. We fix the following: \begin{enumerate}
\item a separable Noether normalization $S := k[t_1, \dots, t_n]_{(t_1, \dots, t_n)} \subseteq R$ such that $S \to R$ is local;
\item $K := k(t_1, \dots, t_n) \subseteq \Frac(R) := L$, the associated separable extension of fraction fields;
\item elements $g_1, \dots, g_r \in R$ such that $R$ is a localization of an $S$-algebra which is generated by $g_1, \dots, g_r$ as an $S$-module;
\item monic polynomials $Q_j \in k[t_1, \dots, t_n][Z]$ such that $Q_j(g_j) = 0$ and $Q_j'(g_j) \neq 0$.
\end{enumerate} Lastly, we denote $d_0 = \prod_j Q_j'(g_j)$.
\end{notation}

Although $D^{[\alpha]}(f) \in k[t_1, \dots, t_n]$ whenever $f \in k[t_1, \dots, t_n]$ and $\alpha \in \N^n$, it is not generally the case that $D^{[\alpha]}(R) \subseteq R$ in the setting above. In what follows, we control these denominators to define operators on $R$ itself.

\begin{proposition} \label{lem: Psi}
Let $d \in R$ be a nonzero element divisible by $d_0$. There is a $k$-algebra homomorphism $$\Psi: R \longrightarrow R\llb T_1, \dots, T_n \rrb, \qquad \Psi(f) = \sum_{\alpha \in \N^n} d^{2|\alpha|}D^{[\alpha]}(f) T^\alpha,$$ extending the map $\Psi_K: K \longrightarrow L\llb T_1, \dots, T_n \rrb$ given by $t_i \mapsto t_i + d^2T_i$. Furthermore, given $f \in R$, this homomorphism satisfies $$\Psi(f) - f \in d(T_1, \dots, T_n)R\llb T_1, \dots, T_n\rrb.$$
\end{proposition}
\begin{proof}
First, we may define an extension $\Psi_L: L \to L\llb T_1, \dots, T_n \rrb$ by composing the $\Theta_L$ of \cref{prop:HSTaylorExpansion} with coordinate-wise multiplication by $d^2$, that is, the substitution $T_i \mapsto d^2T_i$. This ensures $\Psi_L(f) \equiv f \bmod (T_1, \dots, T_n)$ for all $f \in L$. We must show that $\Psi_L(R) \subseteq R{\llb T_1, \dots, T_n \rrb}$ to prove the first assertion. Towards this end, if $f \in k[t_1, \dots, t_n]$, then $\Psi_L(f) = f(t + d^2T) \in R[T_1, \dots, T_n]$, and if $s \in k[t_1, \dots, t_n] \setminus (t_1, \dots, t_n)$, then $s$ is a unit in $R$ since $S \to R$ is local, so $\Psi_L(s) \in s + (T_1, \dots, T_n)R\llb T_1, \dots, T_n \rrb$ is a unit in $R \llb T_1, \dots, T_n \rrb$. Thus, as $\Psi_L$ is a ring homomorphism, it suffices to prove $\Psi_L(g_j) \in R\llb T_1, \dots, T_n \rrb$ for all $j$. 

Towards this end, we fix $g = g_j$, write $Q = Q_j$, and consider $$\tilde{Q}(T, Z) = \frac{Q(t + d^2T, g + dZ)}{dQ'(g)}.$$ Taylor expanding the numerator, the constant term is $Q(t, g) = 0$, the term which is linear in $Z$ and constant in $T$ is $dQ'(g)Z$, and every other term is divisible by $d^2$. Since $d$ is divisible by $Q'(g)$ in $R$, it follows that $$\tilde{Q} \in R[T_1, \dots, T_n][Z] \subseteq R\llb T_1, \dots, T_n \rrb[Z].$$ Furthermore, we may write $$\tilde{Q}(T, Z) = Z + \frac{d}{Q'(g)}G(T, Z),$$ where $G(T, Z)$ has coefficients in $R[T_1, \dots, T_n]$ and lies in the ideal $(T_1, \dots, T_n, Z^2)$. It follows that $\tilde{Q}(T, 0) \in (T_1, \dots, T_n)$ and that $\tilde{Q}_Z(T, 0) \in 1 + (T_1, \dots, T_n)$ is a unit of $R\llb T_1, \dots, T_n \rrb$, so by Hensel's lemma \cite[Theorem 7.3]{eisenbud95} applied to the complete ring $R\llb T_1, \dots, T_n \rrb$, there is a unique $$\psi(T) \in (T_1, \dots, T_n)R\llb T_1, \dots, T_n \rrb$$ such that $Q(t + d^2T, g + d\psi(T)) = 0$. Now, we have $$\Psi_L(Q(t, g)) = Q(t + d^2T, \Psi_L(g)) = 0,$$ so $\Psi_L(g)$ and $g + d\psi(T)$ are both roots of the polynomial $Q(t + d^2T, Z) \in L\llb T_1, \dots, T_n \rrb[Z]$ which are congruent to $g$ modulo $(T_1, \dots, T_n)$. Since $Q'(g) \neq 0$ is a unit in $L$, the uniqueness part of Hensel's lemma implies that $\Psi_L(g) = g + d\psi(T) \in R\llb T_1, \dots, T_n \rrb$, as required.

Since $\Psi(t_i)-t_i=d^2T_i$ and $\Psi(g_j)-g_j=d\psi_j(T)$, the congruence holds on generators, and hence on all of $R$, since $\Psi$ is a ring homomorphism. This concludes the proof.
\end{proof}

\begin{corollary} \label{cor: PsiLocalization}
The map $\Psi$ defined in \cref{lem: Psi} extends to any localization $W^{-1}R$ of $R$, in the sense that $\Psi_L(W^{-1}R) \subseteq W^{-1}R\llb T_1, \dots, T_n \rrb$, and this map still satisfies $$\Psi(f) - f \in d(T_1, \dots, T_n)W^{-1}R\llb T_1, \dots, T_n \rrb$$ for all $f \in W^{-1}R$. Furthermore, if $D$ is an integral Weil divisor on $\Spec R$ supported in $V(d)$, then $\Psi_L(R(D)) \subseteq R(D)\llb T_1, \dots, T_n \rrb$.
\end{corollary}

\begin{proof}
The composite map $$R \longrightarrow R\llb T_1, \dots, T_n \rrb \longrightarrow W^{-1}R\llb T_1, \dots, T_n \rrb$$ sends every $s \in W$ to a unit; it follows from \cref{lem: Psi} that if $s \in W$, then $$\Psi(s) - s \in (T_1, \dots, T_n)W^{-1}R\llb T_1, \dots, T_n \rrb.$$ Thus, $\Psi$ extends uniquely to $W^{-1}R$. Moreover, $$\Psi(r/s) - (r/s) = \frac{s\big(\Psi(r) - r\big) - r\big(\Psi(s) - s\big)}{s\Psi(s)} \in d(T_1, \dots, T_n)W^{-1}R\llb T_1, \dots, T_n \rrb$$ by the congruence part of \cref{lem: Psi}. This concludes the proof of the first assertion.

By the definition of $R(D)$, it suffices to prove the second assertion after localizing at every height-one prime ideal $\mathfrak{p}$ of $R$. Suppose the prime divisor $V(\mathfrak{p}) \subseteq \Spec R$ appears in $D$ with coefficient $a \in \Z$. Then $R(D)_\mathfrak{p} = u^{-a} R_\mathfrak{p}$, where $u$ is a uniformizer of the DVR $R_\mathfrak{p}$. If $a = 0$, then this follows directly from the first assertion. Thus we may assume $a \neq 0$, so that $\mathfrak{p} \in \Supp D \subseteq V(d)$ and $u$ divides $d$ in $R_{\mathfrak{p}}$. The congruence part of the first assertion, applied to $R_{\mathfrak{p}}$, implies that $\Psi(u) = u(1 + v)$ for some $v \in (T_1, \dots, T_n)R_{\mathfrak{p}}\llb T_1, \dots, T_n \rrb$. In particular, $1 + v$ is a unit, so that $$\Psi_L(u^{-a} R_{\mathfrak{p}}) = u^{-a}(1 + v)^{-a}\Psi(R_{\mathfrak{p}}) \subseteq u^{-a} R_{\mathfrak{p}}\llb T_1, \dots, T_n \rrb,$$ as required.
\end{proof}

\subsection*{A Colength Bound} By applying these operators to appropriate splitting ideals in the setting of \cref{not:setup}, we are able to give uniform lower bounds on their colengths. This is the primary positive-characteristic input into the proof of the forward direction of \cref{thm: main}. 
\begin{definition}[{\cite[Definition 3.3]{bst12}}]
Let $\Delta$ be an effective $\Q$-divisor on $\Spec R$. For a positive integer $e$, the \textit{$e$-th splitting ideal} of $(R, \Delta)$ is defined as $$I^\Delta_e = \{ f \in R \;|\; \varphi(F^e_* f) \in \mathfrak{m} \text{ for all } \varphi \in \operatorname{Hom}_R(F^e_*(R(\lceil (p^e - 1)\Delta \rceil)), R) \}.$$ When $\Delta = 0$, we simply write $I_e = I^0_e$.
\end{definition}

Using \cref{cor: PsiLocalization} above, we conclude that splitting ideals are invariant under the differential operators defined by \cref{lem: Psi}. This is similar to \cite[Proposition 5.10]{BJN19}.

\begin{corollary} \label{cor: splittingidealinclusions}
In the setting of \cref{not:setup}, let $\Delta$ be an effective $\Q$-divisor on $\Spec R$, and let $d \in R$ be a nonzero element divisible by $d_0$ with $\Supp \Delta \subseteq V(d)$. If $|\alpha| < p^e$, then $$d^{2|\alpha|}D^{[\alpha]}(I^\Delta_e) \subseteq I^\Delta_e.$$
\end{corollary}
\begin{proof}
    It follows from \cref{cor: PsiLocalization} that $d^{2|\alpha|}D^{[\alpha]}(-)$ defines a map $$ R(\lceil (p^e - 1)\Delta \rceil)  \lre R(\lceil (p^e - 1)\Delta \rceil)$$ which maps $R$ into $R$. In addition, the $F$-pushforward $F^e_*\big(d^{2|\alpha|}D^{[\alpha]}\big)$ is $R$-linear by \cref{rem: froblinearity}. Suppose $\varphi: F^e_* R(\lceil (p^e - 1)\Delta \rceil) \to R$ is $R$-linear. Then the composite map $$F^e_* R(\lceil (p^e - 1)\Delta \rceil) \stackrel{F^e_* (d^{2|\alpha|}D^{[\alpha]})}\longrightarrow F^e_*R(\lceil (p^e - 1)\Delta \rceil) \stackrel{\varphi}\longrightarrow R$$ is an $R$-linear map $ F^e_* R(\lceil (p^e - 1)\Delta \rceil) \to R$ as well. Thus, if $f \in I^\Delta_e$, then $\varphi(F^e_* d^{2|\alpha|} D^{[\alpha]}(f)) \in \mathfrak{m}$ for all $R$-linear $\varphi:  F^e_* R(\lceil (p^e - 1)\Delta \rceil) \to R$. It follows that $d^{2|\alpha|} D^{[\alpha]}(f) \in I^\Delta_e$, as required.
\end{proof}
\begin{proposition} \label{prop: colengthbound}
In the setting of \cref{not:setup}, let $\Delta$ and $d$ be as in \cref{cor: splittingidealinclusions}. If $d^N \notin I^\Delta_e$ for some positive integers $e$ and $N < p^e$, then $$\ell_R(R/I^\Delta_e) \geq \binom{\lfloor N/2 \rfloor+ n}{n}.$$
\end{proposition}
\begin{proof}
We will show that the $k$-vector space of polynomials of degree $\leq M := \lfloor N/2 \rfloor$ injects into $R/I_e^{\Delta}$. Fix a $k$-linear combination $f = \sum_{|\beta| \leq M} f_\beta t^\beta$ and suppose, seeking a contradiction, that $f \in I_e^\Delta$ is nonzero. By hypothesis, $d^N \notin I_e^\Delta$, so in particular $I_e^\Delta$ does not contain any units. Thus, $f$ is not a nonzero constant, and we may choose some nonzero $\alpha \in \N^n$ of maximal degree $|\alpha|$ with $f_\alpha \neq 0$. In this case, since $|\alpha| \leq M < p^e$, it follows from \cref{cor: splittingidealinclusions} and \cref{eq:HSmonomial} that $$d^{2|\alpha|}D^{[\alpha]}(f) = d^{2|\alpha|}f_{\alpha} \in I^\Delta_e.$$ On the other hand, $f_\alpha$ is a unit and $2|\alpha| \leq 2M \leq N$, so $d^N \in I^\Delta_e$, which contradicts the hypothesis. Since $R/\frm = k$, it follows that $\ell_R(R/I^\Delta_e) = \dim_k R/I^\Delta_e \geq \binom{M + n}{n}$.
\end{proof}

\begin{theorem} \label{thm:charpmain}
In the setting of \cref{not:setup}, let $\Delta$ be an effective $\Q$-divisor on $\Spec R$ whose coefficients have denominators prime to $p$, and let $d \in R$ be a nonzero element divisible by $d_0$ with $\Supp \Delta \subseteq V(d)$. If $0< \eta < 1$ is a rational number with denominator prime to $p$ such that $(R, \Delta+ \eta \Div(d))$ is strongly $F$-regular, then the $F$-signature of the pair $(R, \Delta)$ satisfies $$s(R, \Delta) \geq \frac{\eta^n}{2^n n!} > 0.$$
\end{theorem}
\begin{proof}
Choose a positive integer $e_0$ so that for multiples $e$ of $e_0$, the divisor $(q - 1)\Delta$ and the number $(q - 1)\eta$ are integral, where $q = p^e$. In this setting, denote $N_e = (q - 1)\eta$. Replacing $e_0$ by a multiple if necessary, it follows from the strong $F$-regularity of the pair $(R, \Delta+ \eta \Div(d))$ and \cite[Chapter 5, Exercise 3.3]{SSbook} that for all multiples $e$ of $e_0$ there are maps $$\psi: F^e_* R((q - 1)(\Delta+ \eta \Div(d))) \lre R$$ such that $\psi(F^e_*1) = 1$. In this case, viewing $$R((q - 1)(\Delta+ \eta \Div(d))) = d^{-N_e}R((q - 1)\Delta),$$ there is an associated map $\tilde{\psi}$ as below, where the vertical map is multiplication by $d^{-N_e}$.
$$\begin{tikzcd}
F^e_*d^{-N_e}R((q - 1)\Delta) \arrow[r, "\psi"]                      & R \\
F^e_* R((q - 1)\Delta) \arrow[u, "\cdot d^{-N_e}"] \arrow[ru, "\tilde{\psi}"'] &
\end{tikzcd}$$
Since $\tilde{\psi}(F^e_*d^{N_e}) = \psi(F^e_*1) = 1$, it follows that $d^{N_e} \notin I_e^{\Delta}$. Since $N_e = (q - 1)\eta < q$, \cref{prop: colengthbound} implies
$$\ell_R(R/I^\Delta_e) \geq \binom{\lfloor N_e/2 \rfloor + n}{n} \geq \frac{(N_e/2)^n}{n!}.$$
The above inequality holds for all multiples $e$ of $e_0$, which, combined with the fact \cite[Theorem 3.11]{bst12} that the $F$-signature is the limit $s(R, \Delta) = \lim_{e \to \infty} \ell_R(R/I^\Delta_e)/p^{en}$, proves the theorem.
\end{proof}

\subsection*{The Forward Direction} We end this section by explaining why \cref{thm:charpmain} implies \cref{thm: main1}.

\begin{theorem} \label{thm: forwardirection} Let $X_\Lambda \to \Spec \Lambda$ be a model of $x \in X$. If the singularity $x \in X$ is of klt type, then there exist a real number $\vare > 0$ and a nonempty open subset $U \subseteq \Spec \Lambda$ so that $$s(\O_{X_b, x_b}) \geq \vare$$ for all closed points $b \in U$.\end{theorem}

\begin{proof}
First, we may replace $\Lambda$ by a finitely generated $\Lambda$-algebra $\Lambda \subseteq \Lambda' \subseteq \C$. Indeed, by \cite[Theorem 5.6]{yao06}, the $F$-signature only depends on the geometric point over $b$, so by Chevalley's theorem, as in the proof of \cite[Chapter 6, Proposition 1.23]{SSbook}, the existence of a nonempty open subset of $\Spec \Lambda'$ satisfying the conclusion of \cref{thm: forwardirection} or \cref{thm: conversedirection} implies the existence of one for $\Spec \Lambda$. 

We may shrink $X$ to assume there is an effective divisor $\Delta$ on $X = \Spec A$ making $(X, \Delta)$ klt. Fix a finite map $\C[t_1, \ldots, t_n] \subset A$ centered at $x$ so that $A$ is generated by $g_1, \dots, g_r$ as a $\C[t_1, \dots, t_n]$-module. If we denote by $Q_j$ the monic minimal polynomials of $g_j$ over the polynomial ring, we may extend $\Lambda$ to assume that $\Delta$, the $t_i$, the $g_j$ and the $Q_j$ are all defined over $\Lambda$, with $Q_j \in \Lambda[t_1, \dots, t_n][Z]$ monic, $Q_j(g_j) = 0$, $A_\Lambda = \sum_j \Lambda[t_1, \dots, t_n]g_j$, and the $t_i$ vanishing along $x_\Lambda$ \cite[Chapter 6, Section 1.3 and Exercise 1.7]{SSbook}.

   Let $h \in A$ be a nonzero element vanishing on $\Supp \Delta$, and let $$d = h\prod_j Q_j'(g_j) \in A,$$ which is nonzero since each $Q_j$ is separable in characteristic zero, and satisfies $\operatorname{Supp}\Delta\subseteq V(d)$ since $h$ vanishes on $\Supp \Delta$. Let $0 < \eta < 1$ be a rational number so that $(X, \Delta + \eta\Div(d))$ is klt. Such $\eta$ exists by \cite[Corollary 2.35(2)]{kollarmori}. We may extend $\Lambda$ as above to assume $h$ is defined over $\Lambda$. We then choose a nonempty Zariski open set $U \subset \Spec \Lambda$ so that the following hold for closed points $b \in U$: \begin{enumerate}[ref=(\arabic*)]
\item\label{cond: fwdchar} the characteristic of $k(b)$ does not divide the denominators of $\eta$ and of the coefficients of $\Delta$;
   \item\label{cond: fwdnormal} $X_b$ is geometrically integral and normal of dimension $n$;
\item\label{cond: fwdnoether} $k(b)[t_1, \dots, t_n] \to A_b$ is a Noether normalization of $A_b$;
\item\label{cond: fwdnonzero} $Q_j'(g_j) \neq 0$ in $A_b$ for every $j$ and $h_b \neq 0$ in $A_b$, so in particular $d_b \neq 0$ in $A_b$;
\item\label{cond: fwdsfr} $(X_b, \Delta_b + \eta \Div(d_b))$ is strongly $F$-regular at $x_b \in X_b$.
\end{enumerate}

Condition \labelcref{cond: fwdchar} only excludes finitely many characteristics. Condition \labelcref{cond: fwdnormal} holds on a nonempty open set by \cite[Chapter 6, Proposition 1.16(c), (h), (l) and its proof]{SSbook}, and conditions \labelcref{cond: fwdnoether} and \labelcref{cond: fwdnonzero} by \cite[Chapter 6, Lemma 1.12 and Corollary 1.13]{SSbook}, applied to the inclusion $\C[t_1, \dots, t_n] \subseteq A$ and to multiplication by $Q_j'(g_j)$ and $h$ on $A$ (see also \cite[6.9.1]{EGAIV2} and \cite[9.5.6, 9.7.7, 9.9.5]{EGAIV3}). For condition \labelcref{cond: fwdsfr}, we apply \cite[Chapter 6, Theorem 4.23]{SSbook} (see also \cite[Corollary 3.4]{takagi04}) to the klt pair $(X, \Delta + \eta \Div(d))$.

It follows from conditions \labelcref{cond: fwdnoether} and \labelcref{cond: fwdnonzero} that $k(b)[t_1, \dots, t_n] \to A_b$ is a separable Noether normalization. Thus, the local ring $\O_{X_b, x_b}$, the element $d_b$, and the associated Noether normalization are as in \cref{not:setup} and \cref{thm:charpmain}. Moreover, since $h$ vanishes on $\Supp \Delta_\Lambda$, we have $\Supp \Delta_b \subseteq V(d_b)$. By condition \labelcref{cond: fwdchar}, the denominators of $\eta$ and of the coefficients of $\Delta_b$ are prime to the characteristic of $k(b)$, so we may apply \cref{thm:charpmain}:  $$s(\mathcal{O}_{X_b, x_b}) \geq s(\mathcal{O}_{X_b, x_b}, \Delta_b) \geq  \frac{\eta^n}{2^n n!} > 0. $$ Taking $\vare = \frac{\eta^n}{2^n n!}$ thus proves the result.
\end{proof}

\section{Uniform Bounds and klt Type}

In this section, we prove the converse direction of \cref{thm: main}. As in the forward direction, we will first consider the positive characteristic statement, and then conclude the section with the proof of the characteristic zero theorem.

\subsection*{A Cartier Generator} First, we fix a perfect field $k$ of characteristic $p> 0$ and let $(R, \frm)$, $S = k[t_1, \dots, t_n]_{(t_1, \dots, t_n)}$, and $K \subset L$ be as in \cref{not:setup}. To start, we describe an explicit Cartier map $\Phi^e_\omega: F^e_*L \to L$ and relate it to Cartier maps $F^e_*R\to R$ to understand the splitting ideals $I_e$ of $R$.

\begin{notation} \label{not: conversesetup}
We denote $\omega = dt_1 \wedge \cdots \wedge dt_n$ and $K_\omega = \Div_R(\omega)$. Since $\omega$ is a regular $n$-form, $K_\omega$ is effective, so the modules $J_j := R(-jK_\omega) \subseteq R$ are nonzero ideals of $R$. We also fix a nonzero element $c \in J_2 \cap \frm$. For positive integers $e$ and $q = p^e$, since $k$ is perfect and $K \subset L$ is separable, $F^e_* L$ is a free $L$-module with basis $\{F^e_* t^\alpha\}_{\alpha \in \{0, \ldots, q-1\}^n}$, and we denote by $$\Phi^e_{\omega}: F^e_* L \lre L $$ the projection onto the $F^e_* t^{q - 1}$ component, where $t^{q-1} = t_1^{q-1} \cdots t_n^{q-1}$.
\end{notation}

\begin{remark} \label{rem: Cartiertrace}
$\Phi^e_\omega$ is the \textit{Cartier trace relative to the top form $\omega$}, see \cite[Lemma 1.3.6, Proposition 1.3.7, Remark 1.3.12]{BK05}. In particular, if $T^e: F^e_* \Omega^n_{L/k} \to \Omega^n_{L/k}$ is the $e$-th iterate of the Cartier operator, we can explicitly describe \begin{align}\label{eq: ctrace} \Phi^e_\omega(F^e_*f) = T^e(F^e_*(f \omega))/\omega.\end{align} We can define $\Phi^e_\tau$ for other nonzero rational $n$-forms $\tau \in \Omega^n_{L/k}$ in the same way. If $\tau$ is another nonzero rational $n$-form and $\omega = h\tau$ with $h \in L$, then it follows from \cref{eq: ctrace} and the $p^{-e}$-linearity of $\Phi^e_\tau$ that \begin{align} \label{eq: omegatotau} \Phi^e_\omega(F^e_* f) h = \Phi^e_\tau(F^e_* (fh)) \;\;\;\;\;\;\;\;\;\; \Phi^e_\omega(F^e_* f) = \Phi^e_\tau(F^e_* (h^{1 - q} f))\end{align} for $f \in L$.
\end{remark}

\begin{lemma} \label{cartiermaps}
Let $e$ be a positive integer, $q = p^e$, and let $D$ be a Weil divisor on $\Spec R$. Then, the map $a \mapsto \Phi^e_\omega(F^e_*(a \cdot (-)))$ defines an isomorphism of $F^e_*R$-modules $$F^e_*R((1 - q)K_\omega + qD) \stackrel{\sim}{\lre} \Hom_R(F^e_*R, R(D)).$$ In particular, by taking $D = 0$, every $R$-linear map $\psi: F^e_*R \to R$ is of the form $\Phi^e_\omega(F^e_*(a \cdot (-)))$ for a unique $a \in J_{q - 1}$. 
\end{lemma}

\begin{proof}
Since $t_1, \dots, t_n$ is a $p$-basis of $L$, the map $\Phi^e_\omega$ generates $\Hom_L(F^e_*L, L)$ as an $F^e_*L$-module \cite[Chapter 2, Corollary 1.16]{SSbook}, so every $R$-linear map $\psi: F^e_*R \to R(D)$ extends to one of the form $\Phi^e_\omega(F^e_*(a \cdot (-)))$ on $F^e_*L$ for a unique $a \in L$. Thus, we must show that $a \in R((1 - q)K_\omega + qD)$ if and only if $\Phi^e_\omega(F^e_*(aR)) \subseteq R(D)$. By the definition of $R(D)$, it suffices to prove this after localizing at every height-one prime $\mathfrak{p}$. 

Fix such a $\mathfrak{p}$. Since $R_{\frak{p}}$ is a DVR and $k$ is perfect, we can choose a regular $n$-form $\tau$ trivializing $\Omega^n_R \otimes R_{\frak{p}}$, and write $\omega = h\tau$ and $R(D)_\frak{p} = sR_\frak{p}$ with $h, s \in L$. That way, $$R((1-q)K_\omega + qD)_\frak{p} = h^{q-1}s^qR_\frak{p}.$$ By \cref{eq: omegatotau} and the $p^{-e}$-linearity of $\Phi^e_\tau$ (that is, the projection formula), $$\Phi^e_\omega(F^e_*(a \cdot (-))) = s\,\Phi^e_\tau(F^e_*(h^{1-q}s^{-q}a \cdot (-))).$$ Since $\Phi^e_\tau$ generates $\Hom_{R_\frak{p}}(F^e_*R_\frak{p}, R_\frak{p})$ as an $F^e_*R_\frak{p}$-module \cite[Chapter 2, Proposition 1.20]{SSbook}, by uniqueness this map sends $F^e_*R_\frak{p}$ into $sR_\frak{p}$ if and only if $h^{1-q}s^{-q}a \in R_\frak{p}$, as required.
\end{proof}

\begin{remark}
Since $K_\omega$ is a canonical divisor on $R$, \cref{cartiermaps} with $D = 0$ gives an explicit form of the isomorphism $\Hom_R(F^e_*R, R) \cong F^e_* R((1 - q)K_R)$ of \cite[Chapter 3, Theorem 2.25]{SSbook}; both are induced by $T^e$ on the regular locus, see \cref{eq: ctrace}.
\end{remark}

\subsection*{Log Discrepancies} Using the Cartier map constructed above, we relate the exclusion $c^N \notin I_e$ to the positivity of log discrepancies. That way, we will be able to use the non-positivity of certain log discrepancies in the non-klt-type case to obtain inclusions $c^N \in I_e$, which we will use to bound $\ell_R(R/I_e)$ from above in \cref{thm: charpconverse}.

\begin{definition} \label{def: limitinglogdisc}
Given a prime divisor $E$ over $\Spec R$, appearing on a normal birational model $Y \to \Spec R$, we define the \textit{log discrepancy} of $E$ with respect to $\omega$ to be $A_{\omega}(E) = 1 + \operatorname{coeff}_E(\Div_Y\omega)$, and the \textit{$m$-th limiting log discrepancy along $E$} to be $$A_m(E) = A_{\omega}(E) - \frac{1}{m}\ord_E(J_m),$$ where $\ord_E(J_m)$ denotes the minimum of $\ord_E$ on $J_m$. Since $\omega$ is a regular $n$-form on $Y$, we have $A_\omega(E) \geq 1$.
\end{definition}

\begin{remark} \label{rem: dFH}
The same definitions make sense in characteristic $0$, when $\Spec R$ is replaced by a normal complex variety $X = \Spec A$ and $\omega$, $K_\omega$ and $J_m$ are defined using a Noether normalization $\C[t_1, \dots, t_n] \subset A$. In this case, $A_m(E)$ is the $m$-th limiting log discrepancy $a_{m, E}(X, 0)$ of \cite[Definition 7.1]{dFH09}, computed with the canonical divisor $K_X = K_\omega$, see \cite[Definitions 2.6, 3.1, and 7.1]{dFH09}.
\end{remark}

\begin{proposition} \label{prop: discrepancybound}
 In the setting of \cref{not: conversesetup}, let $e$, $m$, and $N$ be positive integers with $m - 1 \leq q = p^e$, and suppose $c^N \notin I_e$. If $E$ is a prime divisor over $\Spec R$ and $N \geq q/m$, then $A_m(E) > 0$.
\end{proposition}

\begin{proof}
Since $c^N \notin I_e$, there is an $R$-linear map $\psi: F^e_*R \to R$ with $\psi(F^e_*c^N) = 1$. Using \cref{cartiermaps}, we write $\psi = \Phi^e_\omega(F^e_*(a \cdot (-)))$ with $a \in J_{q - 1}$, so that $\Phi^e_{\omega}(F^e_*(ac^N)) = 1$. Write $v = \ord_E$, and let $R_E$ be the DVR in $L$ associated to $v$. As in the proof of \cref{cartiermaps}, we fix a regular $n$-form $\tau$ so that $\Phi^e_\tau$ generates $\Hom_{R_E}(F^e_*R_E, R_E)$ as an $F^e_*R_E$-module, and $h \in R_E$ satisfying $\omega = h\tau$ in $\Omega^n_{R_E/k}$. For a nonzero element $f \in L$, we may divide by $q$ to write $v(f) = \ell q + \rho$ with $0 \leq \rho < q$. Thus, if $u$ is a uniformizer of $R_E$, factoring $(u^\ell)^q$ out of $F^e_*f$, we have \begin{align*}
q\,v(\Phi^e_\tau(F^e_*f)) &= q\ell + q\,v(\Phi^e_{\tau}(F^e_*(f/u^{\ell q}))) \\&\geq v(f) - \rho \geq v(f) - (q - 1).
\end{align*}

Consider the element $f = h^{1 -q}ac^N \in L$, so that $\Phi^e_\tau(F^e_*f) = \Phi^e_{\omega}(F^e_*(ac^N)) = 1$ by \cref{eq: omegatotau}. Applying the above inequality to $f$ and solving for $v(a)$ gives \begin{align} \label{eq: aestimate} v(a) \leq (q - 1)A_{\omega}(E) - N v(c), \end{align} using the fact that $v(h) = A_\omega(E) - 1$. On the other hand, we claim the following.

\noindent\textit{Claim.} $\ord_E(J_m) \leq v(c) + \frac{m}{q}v(a)$.

\noindent Given the claim, we combine it with \cref{eq: aestimate} to obtain $$\frac{1}{m}\ord_E(J_m) \leq \left( \frac{1}{m} - \frac{N}{q}\right) v(c) + \frac{q - 1}{q} A_\omega(E).$$ Subtracting the above from $A_\omega(E)$ gives $$A_m(E) \geq \frac{A_\omega(E)}{q} + \left( \frac{N}{q}- \frac{1}{m} \right) v(c).$$ Since $A_{\omega}(E) \geq 1$, $v(c) \geq 0$, and $N \geq q/m$, this implies $A_m(E) > 0$, as required. It remains to prove the claim.

To do so, we write $F^e_*a^m = \sum_\beta b_\beta F^e_* t^\beta$ in the basis of $F^e_*L$ over $L$ from \cref{not: conversesetup}, that is, $a^m = \sum_\beta b_\beta^q t^\beta$ with $\beta \in \{0, \dots, q-1\}^n$. It follows from the explicit form of $\Phi_\omega^e$ given there that we may extract these coefficients as $$b_\beta = \Phi^e_\omega(F^e_*(a^m t^{(q - 1) - \beta})),$$ so that $cb_\beta = \Phi^e_\omega(F^e_*(c^qa^m t^{(q-1)-\beta}))$. Since $c^qa^m \in J_{2q + m(q-1)} \subseteq J_{mq + q - 1}$, where the last inclusion uses $q \geq m - 1$, it follows from \cref{cartiermaps} applied to $D = -mK_\omega$, and evaluation along $F^e_*1$, that $c b_\beta \in J_m$ for all $\beta$. Since $c^qa^m = \sum_\beta (cb_\beta)^q t^\beta$ in $L$ and $v(t^\beta) \geq 0$, we may compute $$q\,v(c) + m\,v(a) = v(c^q a^m)  \geq \min_\beta v((cb_\beta)^q t^\beta) \geq q \ord_E(J_m),$$ as required.
\end{proof}

\begin{theorem}\label{thm: charpconverse}
In the setting of \cref{not: conversesetup}, assume $m$ is a positive integer and $E$ is a prime divisor over $\Spec R$ with $A_m(E) \leq 0$, and let $x_2, \dots, x_n \in \frm$ be such that $c, x_2, \dots, x_n$ is a system of parameters of $R$. Then, for every positive integer $e$ with $q = p^e \geq m - 1$, we have $$\frac{\ell_R(R/I_e)}{q^n} \leq \lambda\left(\frac{1}{m} + \frac{1}{q}\right), \qquad \text{where } \lambda = \ell_R(R/(c, x_2, \dots, x_n)).$$ In particular, the $F$-signature of $R$ satisfies $s(R) \leq \lambda/m$.
\end{theorem}

\begin{proof}
We argue as in the proof of \cite[Theorem 3.3]{BST13}. Let $N = \lceil q/m \rceil$. Since $A_m(E) \leq 0$, it follows from \cref{prop: discrepancybound} that $c^N \in I_e$. On the other hand, since $\frm^{[q]} \subseteq I_e$, we may bound $$\ell_R(R/I_e) \leq \ell_R(R/(c^N, x_2^q, \dots, x_n^q)) \leq N q^{n-1}\lambda \leq \left(\frac{q}{m} + 1\right)q^{n-1}\lambda.$$ The second inequality above is obtained by filtering $R/(c^N, x_2^q, \dots, x_n^q)$ by powers of $c$ and of the $x_i$. This proves the first assertion, and the second follows by letting $e \to \infty$.
\end{proof}

\subsection*{The Converse} Next, using the above, we prove \cref{thm: main2}. Throughout, we fix a closed point $x \in X = \Spec A$, where $A$ is a normal, finite type domain over $\C$.  Let $\Lambda \subseteq \C$ be a finitely generated $\Z$-algebra over which $x$ and $X$ are defined. 

\begin{theorem} \label{thm: conversedirection}
    If $x \in X$ is not of klt type, then for any $\vare > 0$ there is a nonempty open set $U \subseteq \Spec \Lambda$ such that $s(\O_{X_b, x_b}) < \vare$ for all closed points $b \in U$.
\end{theorem}

\begin{proof}
Suppose $x \in X$ is not of klt type, and let $\vare > 0$. We will exhibit a nonempty open subset of $\Spec \Lambda$ such that the reductions of $x \in X$ over this subset have $F$-signature less than $\vare$. We fix a finite map $\C[t_1, \dots, t_n] \subseteq A$ centered at $x$, so that $A$ is generated by $g_1, \dots, g_r$ as a $\C[t_1, \dots, t_n]$-module, with monic minimal polynomials $Q_j$. We fix the following: \begin{enumerate}
    \item $\omega = dt_1 \wedge \cdots \wedge dt_n$ and $K_\omega = \Div(\omega)$;
    \item a nonzero element $c \in A(-2K_\omega)$ vanishing at $x$, and general linear combinations $x_2, \dots, x_n$ of $t_1, \dots, t_n$, so that $A/(c, x_2, \dots, x_n)$ is a $\C$-vector space of finite dimension $\lambda$;
    \item an integer $m \geq 2$ with $m > \lambda/\varepsilon$.
\end{enumerate}

By \cite[Theorem 5.4]{dFH09}, there is an effective $\Q$-divisor $\Delta$ on $X$ with $\lfloor \Delta \rfloor = 0$ and a log resolution $\pi: Y \to X$ of $(X, \O_X(-mK_\omega))$ such that $A_{(X, \Delta)}(E) = A_{m}(E)$ for all prime exceptional divisors $E$ on $Y$. Since $X$ is not of klt type, the pair $(X, \Delta)$ is not klt. It therefore follows from \cite[Corollary 2.32]{kollarmori} that there is some prime divisor $E$ on $Y$ whose image contains $x$, necessarily exceptional since $\lfloor \Delta \rfloor = 0$, for which $A_m(E) \leq 0$.

 As in the proof of \cref{thm: forwardirection}, we may replace $\Lambda$ by a finite-type extension to assume all the above data are defined over $\Lambda$. We then choose a nonempty open subset $U \subseteq \Spec \Lambda$ so that the following hold for closed points $b \in U$:
\begin{enumerate}[ref=(\arabic*)]
    \item\label{cond: convnormal} $X_b$ is geometrically integral, normal, and $n$-dimensional;
    \item\label{cond: convnoether} $k(b)[t_1, \dots, t_n] \subseteq A_b$ is a separable Noether normalization;
    \item\label{cond: convres} $\pi_b: Y_b \to X_b$ is a resolution of $X_b$, and $E_b$ is a prime divisor on $Y_b$ with $x_b \in \pi_b(E_b)$ and $A_m(E_b) = A_m(E) \leq 0$;
    \item\label{cond: convsheaves} $(\O_{X}(-jK_\omega))_b = \O_{X_b}(-jK_{\omega_b})$ for $j = 2, m$;
    \item\label{cond: convlength} $A_b/(c_b, x_2, \dots, x_n)$ is a $k(b)$-vector space of dimension $\lambda$.
\end{enumerate}

We can find an open set $U$ satisfying conditions \labelcref{cond: convnormal} and \labelcref{cond: convnoether} as in the proof of the forward direction, and condition \labelcref{cond: convlength} by generic freeness \cite[Chapter 6, Theorem 1.8]{SSbook}. Resolutions spread out to an open subset of $U$ by \cite[Chapter 6, Theorem 3.2]{SSbook}, so we may assume the first part of \labelcref{cond: convres} holds. To justify condition \labelcref{cond: convsheaves}, we note that since $\pi: Y \to X$ is a proper birational morphism of normal varieties, we have \begin{align} \label{eq: divisorialmode} H^0(Y, \O_Y(-j\widetilde{K}_{\omega})) = A(-jK_{\omega}) \end{align} for all nonnegative integers $j$, where $\widetilde{K}_{\omega}$ is the strict transform of $K_{\omega}$. We apply generic freeness, in the form of \cite[Chapter 6, Theorem 1.8]{SSbook}, to assume that the sheaves $\O_{Y_\Lambda}(-j\widetilde{K}_{\omega_{\Lambda}})$ are flat over $U$ for $j = 2, m$. Thus, by \cite[Chapter 6, Lemma 1.14]{SSbook}, we may shrink $U$ once more to have $$H^0(Y_{\Lambda}, \O_{Y_\Lambda}(-j\widetilde{K}_{\omega_{\Lambda}})) \otimes_\Lambda k(b) = H^0(Y_b,  \O_{Y_b}(-j\widetilde{K}_{\omega_{b}}))$$ for closed points $b \in U$. Since \cref{eq: divisorialmode} also holds for $\pi_b: Y_b \to X_b$, this $U$ satisfies condition \labelcref{cond: convsheaves}. Since $A(-mK_\omega)\O_Y$ is invertible and compatible with base change to $k(b)$ by \labelcref{cond: convsheaves}, and $\Div_{Y_b}(\omega_b) = (\Div_Y \omega)_b$ after shrinking $U$, it follows that $A_m(E_b) = A_m(E) \leq 0$. Thus $U$ satisfies \labelcref{cond: convnormal}--\labelcref{cond: convlength}.

It follows from conditions \labelcref{cond: convnormal}, \labelcref{cond: convnoether}, and \labelcref{cond: convsheaves} above that, for $b \in U$, the local ring $\O_{X_b, x_b}$ and the element $c_b \in \O_{X_b, x_b}$ are in the setting of \cref{not: conversesetup}, and by \labelcref{cond: convlength} that $c_b, x_2, \dots, x_n$ is a system of parameters of $\O_{X_b, x_b}$ with $\ell(\O_{X_b, x_b}/(c_b, x_2, \dots, x_n)) \leq \lambda$. Thus, by \cref{thm: charpconverse}, condition \labelcref{cond: convres}, and the assumption that $m > \lambda/\varepsilon$, we obtain $$s(\O_{X_b, x_b}) \leq \frac{\lambda}{m} < \varepsilon$$ for all closed points $b \in U$. 
\end{proof}

\subsection*{Global Case} Having completed the proof of \cref{thm: main}, we conclude this section by proving the global analogue, \cref{cor: global_version}. Throughout, we fix a normal projective variety $Y$. If $L$ is a line bundle on $Y$, we denote the affine cone by $$C(Y, L) := \Spec \bigoplus_{m \geq 0} H^0(Y, L^m)$$ with cone point $o \in C(Y, L)$.

\begin{lemma} \label{lem: log_fanovsklt_type}
If $A$ is an ample Cartier divisor on $Y$, then the following are equivalent. \footnote{We suspect that this is well-known, but we sketch a short argument for lack of direct reference.}
\begin{enumerate}[ref=(\arabic*)]
\item\label{lem: log_fanovsklt_type1} There exists a boundary $\Delta$ making $(Y, \Delta)$ a klt pair, and $A \sim_\Q -r(K_Y + \Delta)$ for some $r \in \Q_{> 0}$;
\item\label{lem: log_fanovsklt_type2}  The cone point $ o \in C(Y, A) $ is a klt type singularity.
\end{enumerate}  
\end{lemma}
\begin{proof}
\labelcref{lem: log_fanovsklt_type1} implies \labelcref{lem: log_fanovsklt_type2} follows from \cite[Lemma 3.1]{kollarsings}. Conversely, we follow the proof of \cite[Corollary 3]{BGLM24}. Since $C(Y, A)$ is of klt type near $o$ and the klt type locus is a $\G_m$-invariant open subset of $C(Y, A)$, it follows that $C(Y, A)$ is of klt type everywhere. The local class group of $o \in C(Y, A)$ is isomorphic to the global class group of $C(Y, A)$ by \cite[Lemma 2.10]{BM24}, the latter of which is finitely generated by \cite[Proposition 2.12]{BGLM24}. Therefore, it follows from \cite[Proposition 2.14]{BGLM24} that $C(Y, A)$ is a Mori dream space whose Cox ring has klt type spectrum. Applying \cite[Theorem 6]{BGLM24} to the $\G_m$-quotient $p$, we see that $Y$ is also a Mori dream space with klt type Cox spectrum. It follows from the criterion \cite[Corollary 5.3]{GOST15} that $Y$ is of Fano type.

Let $\Delta_0 \geq 0$ be a $\Q$-divisor on $Y$ so that $(Y, \Delta_0)$ is a klt pair and $-K_Y - \Delta_0$ is ample. We choose a large integer $m$ and a small rational number $\vare > 0$ so that $-K_Y- \Delta_0 - \vare A$ is ample, and $m(-K_Y - \Delta_0 - \vare A)$ is very ample. Using \cite[Theorem 4.8]{kollar97singsofpairs}, we may replace $\Delta_0$ by $\Delta = \Delta_0 + \frac{1}{m}G$ for a general member $G \in |m(-K_Y- \Delta_0 - \vare A)|$ to assume $(Y, \Delta)$ is klt. This pair then satisfies \labelcref{lem: log_fanovsklt_type1} with $r = 1/\vare$. 
\end{proof}

\begin{corollary} \label{cor: section_global_version}
Let $Y$ be a normal projective variety over $\C$ and $A$ be an ample Cartier divisor on $Y$. If $Y_\Lambda \to \Spec \Lambda$ is a model of $Y$, then the following hold. 

\begin{enumerate}
\item If $(Y, \Delta)$ is a klt pair and $A \sim_\Q r(-K_Y - \Delta)$ is an ample Cartier divisor on $Y$ for some $r \in \Q_{> 0}$, there is a nonempty open subset $U \subseteq \Spec \Lambda$ and a real number $\varepsilon > 0$ so that $s(\O_{C(Y_b, A_b), o_b}) \geq \varepsilon$ for all $b \in U$.
\item If $Y$ is not of Fano type, $A$ is an ample Cartier divisor on $Y$, and $\varepsilon > 0$ is a real number, then there exists a nonempty open subset $U \subseteq \Spec \Lambda$ so that $s(\O_{C(Y_b, A_b), o_b}) < \varepsilon$ for all $b \in U$.
\end{enumerate}
\end{corollary}
\begin{proof}
 Using \cref{lem: log_fanovsklt_type}, this is immediate from \cref{thm: main} (alternatively, \cref{thm: forwardirection} and \cref{thm: conversedirection}).
\end{proof}

\section{Pure Descent of klt Type Singularities}

Finally, we will prove \cref{thm: puredescent} here. But first let us begin with the following lemma.

\begin{lemma} \label{lem: purefunctionals}
Let $(R,\frm) \to (S,\mathfrak{n})$ be a pure local homomorphism. Given a finite $R$-submodule $M \subseteq S$, there exist elements $a_1,\ldots,a_t \in M$ and forms $\beta_1,\ldots,\beta_t \in \Hom_R(M,R)$ such that
\[
 u \equiv \sum_{j=1}^t a_j \beta_j(u) \pmod{\mathfrak{n}}\qquad \forall u \in M.
\]
For each $u \in M$ that is invertible in $S$, one of the values $\beta_j(u)$ is invertible in $R$.
\end{lemma}

\begin{proof}
Choose the free summand $F$ of $M$ that has largest rank. Write $M=F\oplus N$.
Suppose $u\in N$ is a unit of $S$. The $R$-linear map $R\to S$, $a\mapsto au$, is pure since it is the composite of the given map with multiplication by $u$. It factors through $N$, so $Ru\hookrightarrow N$ is pure as a map of $R$-modules. Since $R$ is Noetherian, $N/Ru$ is finitely presented. By \cite[Corollary 5.2]{HR76}, \cite[\stacktag{058L}]{stacks-project}, $N$ has $Ru$ as a direct summand. But $Ru\cong R$, so $F\oplus Ru$ would be a free summand of larger rank. This proves $N\subseteq\mathfrak{n}$.

Let $a_1,\dots,a_t$ be a basis of $F$. Denote the corresponding coordinate maps on $F$ by $\beta_j$. Extend them to $M$ by defining zero on $N$. Now $u - \sum_j a_j\beta_j(u) \in N \subseteq \mathfrak{n}$ for each $u \in M$. If $\beta_j(u) \in \frm$ for every $j$, then $u \in \frm F + N \subseteq \mathfrak{n}$. This proves the last assertion.
\end{proof}

To use \cref{lem: purefunctionals} in characteristic zero we need a finite $A$-module $M \subseteq B$ before choosing $b$ or $q$. As before, the rings are defined over a finitely generated $\Z$-subalgebra $\Lambda \subset \C$. We enlarge $\Lambda$ as needed and take $b$ closed in a dense open set of $\Spec\Lambda$. Write $p = \operatorname{char} k(b)$ and $q = p^e$. For $r\in\Z_{>0}$ with $r(K_B + \Delta)$ principal, we consider $e$ such that $r| (q - 1)$.

Now we choose $M$. It suffices to bound the poles of $\theta(F^e_*f)$ at infinity for $f \in A_b$. The pullback of $f$ to $Z_b$ is regular. The generators $\theta$ correspond to rational sections of $\O_{Z_b}((1-q)K_{Z_b})$. We choose these sections with poles bounded by $a(q-1)H_b$ for a single integer $a$. Since $a(q - 1) \leq qa$, the formula at height one in \cref{cartiermaps} gives
\[
 \theta(F^e_*A_b) \subseteq H^0(Z_b,\O_{Z_b}(aH_b)) = M_b,\qquad M = H^0(Z,\O_Z(aH)).
\]
Properness gives the finiteness of $M$. For similar use of bounded poles, see \cite[Section 8.4, Lemma 8.4.6]{BP09}. See also \cite[Proposition 3]{Anderson2000}, \cite[Proposition 4.5]{Bli13}, and \cite[Proposition 7.6]{BS19}.

\begin{lemma} \label{lem: finitemodule}
Let $A \subseteq B$ be normal affine $\C$-domains. Suppose $\Delta \geq 0$ is a $\Q$-divisor on $\Spec B$ and $r(K_B + \Delta)$ is principal. Then there exist a finite $A$-submodule $M$ of $B$ and, after spreading out, a dense open set $V \subseteq \Spec\Lambda$ such that for $b \in V$ closed and $q = p^e \equiv 1\pmod r$ one can choose $B_b$-module generators of $\Hom_{B_b}(F^e_*B_b((q - 1)\Delta_b),B_b)$ which carry $F^e_*A_b$ into $M_b$.
\end{lemma}

\begin{proof}
Let $z_1,\ldots,z_s$ be nonzero generators of $B$ as a $\C$-algebra. In the resulting embedding $\Spec B \subseteq \mathbb{P}^s_A$, we have $z_i = X_i/X_0$. Let $Z$ be the normalization of the closure and $H$ the pullback of $(X_0 = 0)$. So $z_i \in H^0(Z,\O_Z(H))$ and $Z \setminus \Supp H = \Spec B$, by normality of $B$. The map $Z \to \Spec A$ is projective since normalization is finite \cite[\stacktag{035S}]{stacks-project}.

Use the same rational top form $\omega$ for $K_B = \Div_B(\omega)$ and $K_Z = \Div_Z(\omega)$. Suppose $g \in \Frac(B)^\times$ has divisor $r(K_B + \Delta)$ on $\Spec B$. Thus $\Div_Z(g) - rK_Z$ is effective off $H$. We may replace $H$ by a positive multiple so that \begin{align} \label{eq:sec4ineq} \Div_Z(g) - rK_Z \geq -rH \end{align} The poles of each $z_i$ are still bounded by $H$. In a monomial with each exponent less than $q$, this gives a pole bound of $s(q-1)H$. We take $M = H^0(Z,\O_Z((s+1)H)) \subseteq B$, where the inclusion is restriction to $\Spec B$. It is a finite $A$-module by \cite[\stacktag{02O6}]{stacks-project}.

On the smooth locus $U = Z_{\rm sm}$, the inequality \cref{eq:sec4ineq} implies that $g\omega^{-r}$ is a section of $\omega_U^{-r}(rH)$. On $U \cap \Spec B$ this section has divisor $r\Delta$. The section spreads out with this equality. We may spread out this section and \cref{eq:sec4ineq}, as well as shrink $\Spec\Lambda$ to assume that $A_b \subseteq B_b$, the fibers $Z_b$ are normal and integral, and $Z_b \setminus U_b$ has codimension at least two \cite[Chapter 6, Lemma 1.12 and Proposition 1.16(c), (h), (l)]{SSbook}. We may take $U_b$ smooth and all $z_{i,b}$ nonzero. The affine chart is still $\Spec B_b$ with effective Cartier boundary $H_b$. We obtain $\Div_{Z_b}(g_b)-rK_{Z_b}\geq-rH_b$ and $\Div_{B_b}(g_b)=r(K_{B_b}+\Delta_b)$, where $K_{Z_b}$ and $K_{B_b}$ are defined by $\omega_b$.

For the model $M_\Lambda=H^0(Z_\Lambda,\O_{Z_\Lambda}((s+1)H_\Lambda))$, generic freeness allows us to apply \cite[Chapter 6, Lemma 1.14]{SSbook} after shrinking the base. We then have $M_\Lambda\otimes_\Lambda k(b)=H^0(Z_b,\O_{Z_b}((s+1)H_b))$ for each closed $b$ in this open set.

For a closed $b$ in this open set and $q\equiv 1\pmod r$, write $q-1=rd$. All Cartier maps in the statement are of the form $F^e_*v\mapsto \Phi^e_{\omega_b}(F^e_*(xv))$, where
\[
 x\in B_b\bigl((1-q)(K_{B_b}+\Delta_b)\bigr)
   =B_b\bigl(-d\,\Div_{B_b}(g_b)\bigr)=g_b^d B_b.
\]
This follows from \cref{cartiermaps} in codimension one. Every $h\in B_b$ is a polynomial in $z_{i,b}$. Since the coefficients of this polynomial have $q$-th roots in $k(b)$, we have $h=\sum_\alpha h_\alpha^q z_b^\alpha$ for some $h_\alpha\in B_b$, with $0\le\alpha_i<q$. So $g_b^dh=\sum_\alpha h_\alpha^q(g_b^dz_b^\alpha)$. The maps with $x=g_b^d z_b^\alpha$ generate over $B_b$ by linearity.

As $Z_b$ lies over $\Spec A_b$, every $f\in A_b$ is regular on $Z_b$. For $f\ne0$ and $x=g_b^d z_b^\alpha$ we therefore have
\begin{align*}
 \Div_{Z_b}(xf)+(1-q)K_{Z_b}
 &\geq d\,(\Div_{Z_b}(g_b)-r K_{Z_b})-|\alpha|H_b\\
 &\geq -(q-1)H_b-|\alpha|H_b
 \geq -(s+1)(q-1)H_b.
\end{align*}
In particular $xf$ belongs to $\O_{Z_b}((1-q)K_{Z_b}+q(s+1)H_b)$. Its trace belongs to $\O_{Z_b}((s+1)H_b)$ at height one by \cref{cartiermaps}. Thus $\Phi^e_{\omega_b}(F^e_*(xf))\in M_b$.
\end{proof}

\begin{remark}\label{rem:localization}
Let $W\subseteq A_b$ be a multiplicative subset. For each generator $\theta$ in \cref{lem: finitemodule}, the localization $\theta_W$ sends $F^e_*(W^{-1}A_b)$ into $W^{-1}M_b$.  Take $f\in A_b$ and $t\in W$. Multiplying $F^e_*(f/t)$ by $t$ gives $F^e_*(t^{q-1}f)$. Hence $t\,\theta_W(F^e_*(f/t))=\theta(F^e_*(t^{q-1}f))\in M_b$ by the lemma.
\end{remark}

In the notation of \cref{lem: finitemodule}, take a prime $\mathfrak q$ of $B$ over a prime $\mathfrak p$ of $A$. The assumption from now on is that $R=A_{\mathfrak p}\to B_{\mathfrak q}=S$ is pure. We write $R_b$ for $(A_b)_{\mathfrak p_b}$ and $S_b$ for $(B_b)_{\mathfrak q_b}$. Fix an effective $\Q$-divisor $\Delta$ on $\Spec S$ and a positive integer $r$ for which $r(K_S+\Delta)$ is Cartier.

\begin{proposition} \label{prop: splittingcomparison}
After enlarging $\Lambda$ and shrinking $\Spec\Lambda$, one has $I_e(R_b)\subseteq I_e^{\Delta_b}(S_b)\cap R_b$ for every closed $b$ and every $q=p^e$ with $r|(q-1)$.
\end{proposition}

\begin{proof}
Extend $\Delta$ to $\Spec B$. The divisor $r(K_S + \Delta)$ is principal since $S$ is local. We may therefore invert an element of $B \setminus \mathfrak{q}$ and assume that $r(K_B + \Delta)$ is principal. Take $M \subseteq B$ as in \cref{lem: finitemodule}.

By \cref{lem: purefunctionals}, the residue map $M_{\mathfrak{p}} \to \kappa(\mathfrak{q})$ factors through $R^n$ for some $n$. Write $\beta=(\beta_j)_j$ for the first map. The second map $\gamma$ takes the $j$-th basis vector to the residue of $a_j\in M_{\mathfrak{p}}$. Since $M$ is finitely presented, $\beta$ is defined on $M_w$ over some $A_w$ with $w \notin \mathfrak{p}$ \cite[\stacktag{0583}]{stacks-project}. Choose $w$ so that the elements $a_j$ also belong to $M_w$. Then $\gamma(A_w^n) \subseteq B_w/\mathfrak{q}B_w$. For $u \in M_w$ we have $\gamma\beta(u) = \bar u$ in this ring, as we may check in $\kappa(\mathfrak{q})$.

Fix generators of $M_w$ and spread out the equality on these generators. We may take $w_b \notin \mathfrak{p}_b$ and the element inverted in $B$ outside $\mathfrak{q}_b$. We may also take $A_b/\mathfrak{p}_b \subseteq B_b/\mathfrak{q}_b$ an inclusion of domains \cite[Chapter 6, Lemma 1.12 and Proposition 1.16(h)]{SSbook}.
Localizing gives maps $\beta_b:N \to R_b^n$ and $\gamma_b:R_b^n \to \kappa(\mathfrak{q}_b)$ for $N = (M_b)_{\mathfrak{p}_b}$. For each $u \in N$ the composite is its residue $\bar u$.

Suppose $h \in R_b$ is not in $I_e^{\Delta_b}(S_b)$. Then $\varphi(F^e_*h)$ is a unit for some Cartier map $\varphi$ on the pair. Since the Frobenius module is finitely presented, the generators $\theta$ in \cref{lem: finitemodule} generate also over $S_b$ \cite[\stacktag{0583}]{stacks-project}. At least one of them takes $F^e_*h$ to a unit $u$. By \cref{rem:localization}, $\theta(F^e_*R_b) \subseteq N$. Now $\gamma_b\beta_b(u) = \bar u \ne 0$. So some $\beta_{j,b}(u)$ is a unit in $R_b$. In particular, the $R_b$-linear composite $\beta_{j,b} \circ \theta:F^e_*R_b \to R_b$ takes $F^e_*h$ to a unit. Thus $h \notin I_e(R_b)$.
\end{proof}

Now we are ready to prove \cref{thm: puredescent}. 

\begin{theorem}
    Let $f: \Spec B \to \Spec A$ be a pure morphism of normal affine complex varieties. If $\Spec B$ is of klt type, then so is $\Spec A$.
\end{theorem}

\begin{proof}
Fix a closed point $x\in\Spec A$ with maximal ideal $\mathfrak p$, and set $R=A_{\mathfrak p}$. The ring map $A\to B$ induced by $f$ is pure. So its localization $R\to (A\setminus\mathfrak p)^{-1}B$ is pure. By \cite[Lemma 2.4]{hashimoto23}, there is a prime $\mathfrak q$ of $B$ over $\mathfrak p$ such that $R\to B_{\mathfrak q}$ is pure. Set $S=B_{\mathfrak q}$.

Suppose $R$ is not of klt type. Take a Noether normalization centered at $x$ and let $J_j = A(-jK_\omega)$ be the ideals of \cref{rem: dFH}. Choose $0 \ne c \in J_2 \cap \mathfrak{p}$ and an effective $\Delta$ making $(S,\Delta)$ klt. Choose $m \geq 3$ so that $(S,\Delta + \tfrac{2}{m}\Div(c))$ is klt \cite[Corollary 2.35(2)]{kollarmori}. Take a multiple $r$ of $m$ for which $r(K_S + \Delta)$ is Cartier. For this $m$ there is a prime divisor $E$ over $\Spec R$ with $A_m(E) \leq 0$ \cite[Definition 7.1 and Proposition 7.2]{dFH09}.

Take a log resolution of $J_2$ and $J_m$ carrying $E$. Spread it out with a point of $E$ over $x$, as in the proof of the converse of \cref{thm: main}. We may shrink the base so that the Noether normalization remains separable and $0\ne c_b\in R_b(-2K_{\omega_b})\cap\mathfrak{p}_bR_b$. Then $E_b$ is prime over $\Spec R_b$ with $A_m(E_b)=A_m(E)$ after a further shrinking. On the same open set we may assume that $p \nmid r$, \cref{prop: splittingcomparison} applies, and $(S_b,\Delta_b + \tfrac{2}{m}\Div(c_b))$ is strongly $F$-regular \cite[Corollary 3.4]{takagi04}.

Fix such a $b$. Since $p \nmid r$, we can iterate a splitting of this pair until $q = p^e \equiv 1\pmod r$ \cite[Chapter 5, Proposition 2.2]{SSbook}. Put $N = 2(q - 1)/m$. The domain of the splitting is the Frobenius pushforward of $c_b^{-N}S_b((q-1)\Delta_b)$. Multiplication by $c_b^{-N}$ followed by the splitting gives a map $F^e_*S_b((q - 1)\Delta_b) \to S_b$ taking $F^e_*c_b^N$ to $1$. Thus $c_b^N \notin I_e^{\Delta_b}(S_b)$, so $c_b^N \notin I_e(R_b)$ by \cref{prop: splittingcomparison}. Since $m | q - 1$, we have $q \geq m + 1$ and $N \geq q/m$. \cref{prop: discrepancybound} gives $A_m(E_b) > 0$, a contradiction.

It follows that $A_{\mathfrak{p}}$ is of klt type at every closed point of $\Spec A$. Hence $\Spec A$ is of klt type \cite[Proposition 7.2 and Corollary 7.6]{dFH09}.
\end{proof}

\bibliographystyle{amsalpha}
\bibliography{References}

@article{EGAIV,
  author  = {Grothendieck, Alexander},
  title   = {\'El\'ements de g\'eom\'etrie alg\'ebrique : {IV}.
             \'Etude locale des sch\'emas et des morphismes de sch\'emas,
             quatri\`eme partie},
  journal = {Publications Math\'ematiques de l'IH\'ES},
  volume  = {32},
  year    = {1967},
  pages   = {5--361},
  doi     = {10.1007/BF02732123},
  url     = {https://www.numdam.org/item/PMIHES_1967__32__5_0/}
}

@misc{stacks-project,
  author       = {The {Stacks project authors}},
  title        = {The Stacks project},
  howpublished = {\url{https://stacks.math.columbia.edu}},
  year         = {2026},
}

@misc{pande2025frobeniusversiontiansalphainvariant,
      title={A Frobenius Version of Tians Alpha-Invariant}, 
      author={Suchitra Pande},
      year={2025},
      eprint={2311.00989},
      archivePrefix={arXiv},
      primaryClass={math.AG},
      url={https://arxiv.org/abs/2311.00989}, 
}

@book{kollarsings,
  author = {Koll{\'a}r, J{\'a}nos},
  title = {Singularities of the Minimal Model Program},
  series = {Cambridge Tracts in Mathematics},
  volume = {200},
  publisher = {Cambridge University Press},
  year = {2013}
}

@article {crst18,
    AUTHOR = {Carvajal-Rojas, Javier and Schwede, Karl and Tucker, Kevin},
     TITLE = {Fundamental groups of {$F$}-regular singularities via
              {$F$}-signature},
   JOURNAL = {Ann. Sci. \'Ec. Norm. Sup\'er. (4)},
  FJOURNAL = {Annales Scientifiques de l'\'Ecole Normale Sup\'erieure.
              Quatri\`eme S\'erie},
    VOLUME = {51},
      YEAR = {2018},
    NUMBER = {4},
     PAGES = {993--1016},
      ISSN = {0012-9593,1873-2151},
   MRCLASS = {14B05 (14F35 14G17)},
  MRNUMBER = {3861567},
MRREVIEWER = {Daniel\ Smolkin},
       DOI = {10.24033/asens.2370},
       URL = {https://doi.org/10.24033/asens.2370},
}

@article {bst12,
    AUTHOR = {Blickle, Manuel and Schwede, Karl and Tucker, Kevin},
     TITLE = {{$F$}-signature of pairs and the asymptotic behavior of
              {F}robenius splittings},
   JOURNAL = {Adv. Math.},
  FJOURNAL = {Advances in Mathematics},
    VOLUME = {231},
      YEAR = {2012},
    NUMBER = {6},
     PAGES = {3232--3258},
      ISSN = {0001-8708,1090-2082},
   MRCLASS = {13A35 (13D40 14B05)},
  MRNUMBER = {2980498},
       DOI = {10.1016/j.aim.2012.09.007},
       URL = {https://doi.org/10.1016/j.aim.2012.09.007},
}

@misc{lp26,
      title={On positivity of the limit {F}-signature}, 
      author={Yuchen Liu and Suchitra Pande},
      year={2026},
      eprint={2605.16636},
      archivePrefix={arXiv},
      primaryClass={math.AG},
      url={https://arxiv.org/abs/2605.16636}, 
}

@article {zhuang24,
    AUTHOR = {Zhuang, Ziquan},
     TITLE = {Direct summands of klt singularities},
   JOURNAL = {Invent. Math.},
  FJOURNAL = {Inventiones Mathematicae},
    VOLUME = {237},
      YEAR = {2024},
    NUMBER = {3},
     PAGES = {1683--1695},
      ISSN = {0020-9910,1432-1297},
       DOI = {10.1007/s00222-024-01281-1},
       URL = {https://doi.org/10.1007/s00222-024-01281-1},
}

@article {hw02,
    AUTHOR = {Hara, Nobuo and Watanabe, Kei-Ichi},
     TITLE = {F-regular and {F}-pure rings vs. log terminal and log
              canonical singularities},
   JOURNAL = {J. Algebraic Geom.},
  FJOURNAL = {Journal of Algebraic Geometry},
    VOLUME = {11},
      YEAR = {2002},
    NUMBER = {2},
     PAGES = {363--392},
      ISSN = {1056-3911,1534-7486},
       DOI = {10.1090/S1056-3911-01-00306-X},
       URL = {https://doi.org/10.1090/S1056-3911-01-00306-X},
}

@article {takagi04,
    AUTHOR = {Takagi, Shunsuke},
     TITLE = {An interpretation of multiplier ideals via tight closure},
   JOURNAL = {J. Algebraic Geom.},
  FJOURNAL = {Journal of Algebraic Geometry},
    VOLUME = {13},
      YEAR = {2004},
    NUMBER = {2},
     PAGES = {393--415},
      ISSN = {1056-3911,1534-7486},
       DOI = {10.1090/S1056-3911-03-00366-7},
       URL = {https://doi.org/10.1090/S1056-3911-03-00366-7},
}

@article {tucker12,
    AUTHOR = {Tucker, Kevin},
     TITLE = {{$F$}-signature exists},
   JOURNAL = {Invent. Math.},
  FJOURNAL = {Inventiones Mathematicae},
    VOLUME = {190},
      YEAR = {2012},
    NUMBER = {3},
     PAGES = {743--765},
      ISSN = {0020-9910,1432-1297},
       DOI = {10.1007/s00222-012-0389-0},
       URL = {https://doi.org/10.1007/s00222-012-0389-0},
}

@article {HL02,
    AUTHOR = {Huneke, Craig and Leuschke, Graham J.},
     TITLE = {Two theorems about maximal {C}ohen-{M}acaulay modules},
   JOURNAL = {Math. Ann.},
  FJOURNAL = {Mathematische Annalen},
    VOLUME = {324},
      YEAR = {2002},
    NUMBER = {2},
     PAGES = {391--404},
      ISSN = {0025-5831,1432-1807},
   MRCLASS = {13C14 (13A35 13D07 13H10)},
  MRNUMBER = {1933863},
MRREVIEWER = {\=I.\ \=I.\ Burban},
       DOI = {10.1007/s00208-002-0343-3},
       URL = {https://doi.org/10.1007/s00208-002-0343-3},
}

@misc{xu2026openproblemskstabilityfano,
      title={Open problems in {K}-stability of Fano varieties}, 
      author={Chenyang Xu and Ziquan Zhuang},
      year={2026},
      eprint={2601.15576},
      archivePrefix={arXiv},
      primaryClass={math.AG},
      url={https://arxiv.org/abs/2601.15576}, 
}

@misc{rahulmatthewbertini,
      title={Bertini theorems for {F-}singularities over finite fields}, 
      author={Rahul Ajit and Matthew Bertucci},
      year={2026},
      eprint={},
      archivePrefix={},
      primaryClass={math.AG},
      url={},
      note={private draft}
}

@book {eisenbud95,
    AUTHOR = {Eisenbud, David},
     TITLE = {Commutative algebra},
    SERIES = {Graduate Texts in Mathematics},
    VOLUME = {150},
      NOTE = {With a view toward algebraic geometry},
 PUBLISHER = {Springer-Verlag, New York},
      YEAR = {1995},
     PAGES = {xvi+785},
      ISBN = {0-387-94268-8; 0-387-94269-6},
       DOI = {10.1007/978-1-4612-5350-1},
       URL = {https://doi.org/10.1007/978-1-4612-5350-1},
}

@book {BP09,
    AUTHOR = {B\"ockle, Gebhard and Pink, Richard},
     TITLE = {Cohomological theory of crystals over function fields},
    SERIES = {EMS Tracts in Mathematics},
    VOLUME = {9},
 PUBLISHER = {European Mathematical Society (EMS), Z\"urich},
      YEAR = {2009},
     PAGES = {viii+187},
      ISBN = {978-3-03719-074-6},
   MRCLASS = {14F05 (11G09 14F43 14G10 14G15)},
  MRNUMBER = {2561048},
MRREVIEWER = {Martin\ C.\ Olsson},
       DOI = {10.4171/074},
       URL = {https://doi-org.ezproxy.lib.utah.edu/10.4171/074},
}

@article {Anderson2000,
    AUTHOR = {Anderson, G. W.},
     TITLE = {An elementary approach to {$L$}-functions mod {$p$}},
   JOURNAL = {J. Number Theory},
  FJOURNAL = {Journal of Number Theory},
    VOLUME = {80},
      YEAR = {2000},
    NUMBER = {2},
     PAGES = {291--303},
      ISSN = {0022-314X,1096-1658},
   MRCLASS = {11M38},
  MRNUMBER = {1740516},
MRREVIEWER = {Yoshinori\ Hamahata},
       DOI = {10.1006/jnth.1999.2452},
       URL = {https://doi-org.ezproxy.lib.utah.edu/10.1006/jnth.1999.2452},
}

@article {BS19,
    AUTHOR = {Blickle, Manuel and St\"abler, Axel},
     TITLE = {Functorial test modules},
   JOURNAL = {J. Pure Appl. Algebra},
  FJOURNAL = {Journal of Pure and Applied Algebra},
    VOLUME = {223},
      YEAR = {2019},
    NUMBER = {4},
     PAGES = {1766--1800},
      ISSN = {0022-4049,1873-1376},
   MRCLASS = {13A35 (14B05 14F10)},
  MRNUMBER = {3906525},
MRREVIEWER = {Geoffrey\ D.\ Dietz},
       DOI = {10.1016/j.jpaa.2018.07.007},
       URL = {https://doi-org.ezproxy.lib.utah.edu/10.1016/j.jpaa.2018.07.007},
}

@book {kollarmori,
    AUTHOR = {Koll\'ar, J\'anos and Mori, Shigefumi},
     TITLE = {Birational geometry of algebraic varieties},
    SERIES = {Cambridge Tracts in Mathematics},
    VOLUME = {134},
      NOTE = {With the collaboration of C. H. Clemens and A. Corti,
              Translated from the 1998 Japanese original},
 PUBLISHER = {Cambridge University Press, Cambridge},
      YEAR = {1998},
     PAGES = {viii+254},
      ISBN = {0-521-63277-3},
       DOI = {10.1017/CBO9780511662560},
       URL = {https://doi.org/10.1017/CBO9780511662560},
}

@article{yao06,
    AUTHOR = {Yao, Yongwei},
     TITLE = {Observations on the {$F$}-signature of local rings of characteristic {$p$}},
   JOURNAL = {J. Algebra},
  FJOURNAL = {Journal of Algebra},
    VOLUME = {299},
      YEAR = {2006},
    NUMBER = {1},
     PAGES = {198--218},
       URL = {https://www.sciencedirect.com/science/article/pii/S0021869305004552},
}

@article{EGAIV2,
    AUTHOR = {Grothendieck, Alexander},
     TITLE = {\'{E}l\'{e}ments de g\'{e}om\'{e}trie alg\'{e}brique. {IV}. \'{E}tude locale des sch\'{e}mas et des morphismes de sch\'{e}mas. {II}},
   JOURNAL = {Inst. Hautes \'{E}tudes Sci. Publ. Math.},
    NUMBER = {24},
      YEAR = {1965},
     PAGES = {231},
}

@article{EGAIV3,
    AUTHOR = {Grothendieck, Alexander},
     TITLE = {\'{E}l\'{e}ments de g\'{e}om\'{e}trie alg\'{e}brique. {IV}. \'{E}tude locale des sch\'{e}mas et des morphismes de sch\'{e}mas. {III}},
   JOURNAL = {Inst. Hautes \'{E}tudes Sci. Publ. Math.},
    NUMBER = {28},
      YEAR = {1966},
     PAGES = {255},
}

@book{BK05,
  author = {Brion, Michel and Kumar, Shrawan},
  title = {Frobenius Splitting Methods in Geometry and Representation Theory},
  series = {Progress in Mathematics},
  volume = {231},
  publisher = {Birkh{\"a}user Boston},
  address = {Boston, MA},
  year = {2005},
  doi = {10.1007/b137486}
}

@article {dFH09,
    AUTHOR = {de Fernex, Tommaso and Hacon, Christopher D.},
     TITLE = {Singularities on normal varieties},
   JOURNAL = {Compos. Math.},
  FJOURNAL = {Compositio Mathematica},
    VOLUME = {145},
      YEAR = {2009},
    NUMBER = {2},
     PAGES = {393--414},
      ISSN = {0010-437X,1570-5846},
       DOI = {10.1112/S0010437X09003996},
       URL = {https://doi.org/10.1112/S0010437X09003996},
}

@article{HR76,
    AUTHOR = {Hochster, Melvin and Roberts, Joel L.},
     TITLE = {The purity of the {F}robenius and local cohomology},
   JOURNAL = {Advances in Math.},
    VOLUME = {21},
      YEAR = {1976},
    NUMBER = {2},
     PAGES = {117--172},
}

@article{xu14,
    AUTHOR = {Xu, Chenyang},
     TITLE = {Finiteness of algebraic fundamental groups},
   JOURNAL = {Compos. Math.},
    VOLUME = {150},
      YEAR = {2014},
    NUMBER = {3},
     PAGES = {409--414},
}

@unpublished{SSbook,
    AUTHOR = {Schwede, Karl and Smith, Karen E.},
     TITLE = {Singularities defined by the {F}robenius map},
      NOTE = {Book draft, compiled August 15, 2026. Available at \url{https://github.com/kschwede/FrobeniusSingularitiesBook}},
      YEAR = {2026},
}

@misc{hashimoto23,
    AUTHOR = {Hashimoto, Mitsuyasu},
     TITLE = {Hochster--{E}agon type theorem for {S}erre's {$(S_n)$} condition},
      YEAR = {2023},
      NOTE = {arXiv:2306.14366},
}

@misc{BGO17,
      title={Finiteness of \'etale fundamental groups by reduction modulo $p$}, 
      author={Bhargav Bhatt and Ofer Gabber and Martin Olsson},
      year={2017},
      eprint={1705.07303},
      archivePrefix={arXiv},
      primaryClass={math.AG},
      url={https://arxiv.org/abs/1705.07303}, 
}

@article{TY25,
  title={On the uniform positivity of {$F$}-signature
         under reduction modulo {$p$}},
  author={Shunsuke Takagi and Tatsuki Yamaguchi},
  journal={Journal f{\"u}r die reine und angewandte Mathematik
           (Crelles Journal)},
  year={2026},
  note={To appear},
  doi={10.1515/crelle-2026-0076},
  eprint={2507.16566},
  archivePrefix={arXiv},
  primaryClass={math.AG},
  url={https://arxiv.org/abs/2507.16566}
}

@article {Schoutens05,
    AUTHOR = {Schoutens, Hans},
     TITLE = {Log-terminal singularities and vanishing theorems via
              non-standard tight closure},
   JOURNAL = {J. Algebraic Geom.},
  FJOURNAL = {Journal of Algebraic Geometry},
    VOLUME = {14},
      YEAR = {2005},
    NUMBER = {2},
     PAGES = {357--390},
      ISSN = {1056-3911,1534-7486},
   MRCLASS = {13A35 (14E05)},
  MRNUMBER = {2123234},
MRREVIEWER = {Karen\ E.\ Smith},
       DOI = {10.1090/S1056-3911-04-00395-9},
       URL = {https://doi-org.ezproxy.lib.utah.edu/10.1090/S1056-3911-04-00395-9},
}

@article {BJN19,
    AUTHOR = {Brenner, Holger and Jeffries, Jack and N\'u\~nez-Betancourt,
              Luis},
     TITLE = {Quantifying singularities with differential operators},
   JOURNAL = {Adv. Math.},
  FJOURNAL = {Advances in Mathematics},
    VOLUME = {358},
      YEAR = {2019},
     PAGES = {106843, 89},
      ISSN = {0001-8708,1090-2082},
   MRCLASS = {14F10 (13A35 16S32)},
  MRNUMBER = {4020453},
MRREVIEWER = {Ana\ Bravo},
       DOI = {10.1016/j.aim.2019.106843},
       URL = {https://doi-org.ezproxy.lib.utah.edu/10.1016/j.aim.2019.106843},
}

@article {AL03,
    AUTHOR = {Aberbach, Ian M. and Leuschke, Graham J.},
     TITLE = {The {$F$}-signature and strong {$F$}-regularity},
   JOURNAL = {Math. Res. Lett.},
  FJOURNAL = {Mathematical Research Letters},
    VOLUME = {10},
      YEAR = {2003},
    NUMBER = {1},
     PAGES = {51--56},
      ISSN = {1073-2780},
   MRCLASS = {13A35},
  MRNUMBER = {1960123},
MRREVIEWER = {Karen\ E.\ Smith},
       DOI = {10.4310/MRL.2003.v10.n1.a6},
       URL = {https://doi-org.ezproxy.lib.utah.edu/10.4310/MRL.2003.v10.n1.a6},
}

@article {JS22,
    AUTHOR = {Jeffries, Jack and Smirnov, Ilya},
     TITLE = {A transformation rule for natural multiplicities},
   JOURNAL = {Int. Math. Res. Not. IMRN},
  FJOURNAL = {International Mathematics Research Notices. IMRN},
      YEAR = {2022},
    NUMBER = {2},
     PAGES = {999--1015},
      ISSN = {1073-7928,1687-0247},
   MRCLASS = {13F35 (13H15)},
  MRNUMBER = {4368877},
MRREVIEWER = {\c Serban\ B\u arc\u anescu},
       DOI = {10.1093/imrn/rnaa083},
       URL = {https://doi-org.ezproxy.lib.utah.edu/10.1093/imrn/rnaa083},
}

@article {Bli13,
    AUTHOR = {Blickle, Manuel},
     TITLE = {Test ideals via algebras of {$p^{-e}$}-linear maps},
   JOURNAL = {J. Algebraic Geom.},
  FJOURNAL = {Journal of Algebraic Geometry},
    VOLUME = {22},
      YEAR = {2013},
    NUMBER = {1},
     PAGES = {49--83},
      ISSN = {1056-3911,1534-7486},
   MRCLASS = {13A35 (14F18)},
  MRNUMBER = {2993047},
MRREVIEWER = {Karen\ E.\ Smith},
       DOI = {10.1090/S1056-3911-2012-00576-1},
       URL = {https://doi-org.ezproxy.lib.utah.edu/10.1090/S1056-3911-2012-00576-1},
}

@article {BST13,
    AUTHOR = {Blickle, Manuel and Schwede, Karl and Tucker, Kevin},
     TITLE = {{$F$}-signature of pairs: continuity, {$p$}-fractals and
              minimal log discrepancies},
   JOURNAL = {J. Lond. Math. Soc. (2)},
  FJOURNAL = {Journal of the London Mathematical Society. Second Series},
    VOLUME = {87},
      YEAR = {2013},
    NUMBER = {3},
     PAGES = {802--818},
      ISSN = {0024-6107,1469-7750},
   MRCLASS = {13A35 (13D40 14B05 14F18)},
  MRNUMBER = {3073677},
MRREVIEWER = {Adela\ N.\ Vraciu},
       DOI = {10.1112/jlms/jds070},
       URL = {https://doi-org.ezproxy.lib.utah.edu/10.1112/jlms/jds070},
}

@article {Traves00,
    AUTHOR = {Traves, William N.},
     TITLE = {Tight closure and differential simplicity},
   JOURNAL = {J. Algebra},
  FJOURNAL = {Journal of Algebra},
    VOLUME = {228},
      YEAR = {2000},
    NUMBER = {2},
     PAGES = {457--476},
      ISSN = {0021-8693,1090-266X},
   MRCLASS = {13N15 (13A35 13B40)},
  MRNUMBER = {1764573},
MRREVIEWER = {Karen\ E.\ Smith},
       DOI = {10.1006/jabr.1999.8271},
       URL = {https://doi-org.ezproxy.lib.utah.edu/10.1006/jabr.1999.8271},
}

@article {xz21,
    AUTHOR = {Xu, Chenyang and Zhuang, Ziquan},
     TITLE = {Uniqueness of the minimizer of the normalized volume function},
   JOURNAL = {Camb. J. Math.},
  FJOURNAL = {Cambridge Journal of Mathematics},
    VOLUME = {9},
      YEAR = {2021},
    NUMBER = {1},
     PAGES = {149--176},
      ISSN = {2168-0930,2168-0949},
   MRCLASS = {14B05 (13A18 14E30)},
  MRNUMBER = {4325260},
MRREVIEWER = {Yuchen\ Liu},
       DOI = {10.4310/CJM.2021.v9.n1.a2},
       URL = {https://doi.org/10.4310/CJM.2021.v9.n1.a2},
}

@article{BM24,
    AUTHOR = {Braun, Lukas and Moraga, Joaqu\'{\i}n},
     TITLE = {Iteration of {C}ox rings of klt singularities},
   JOURNAL = {J. Topol.},
    VOLUME = {17},
      YEAR = {2024},
    NUMBER = {1},
       DOI = {10.1112/topo.12321},
       URL = {https://doi.org/10.1112/topo.12321},
}

@article{BGLM24,
    AUTHOR = {Braun, Lukas and Greb, Daniel and Langlois, Kevin and Moraga, Joaqu\'{\i}n},
     TITLE = {Reductive quotients of klt singularities},
   JOURNAL = {Invent. Math.},
    VOLUME = {237},
      YEAR = {2024},
     PAGES = {1643--1682},
       DOI = {10.1007/s00222-024-01280-2},
       URL = {https://doi.org/10.1007/s00222-024-01280-2},
}

@article{GOST15,
    AUTHOR = {Gongyo, Yoshinori and Okawa, Shinnosuke and Sannai, Akiyoshi and Takagi, Shunsuke},
     TITLE = {Characterization of varieties of {F}ano type via singularities of {C}ox rings},
   JOURNAL = {J. Algebraic Geom.},
    VOLUME = {24},
      YEAR = {2015},
    NUMBER = {1},
     PAGES = {159--182},
       DOI = {10.1090/S1056-3911-2014-00641-X},
       URL = {https://doi.org/10.1090/S1056-3911-2014-00641-X},
}

@article {HJSP23,
    AUTHOR = {Hochster, Melvin and Jeffries, Jack and Pandey, Vaibhav and
              Singh, Anurag K.},
     TITLE = {When are the natural embeddings of classical invariant rings
              pure?},
   JOURNAL = {Forum Math. Sigma},
  FJOURNAL = {Forum of Mathematics. Sigma},
    VOLUME = {11},
      YEAR = {2023},
     PAGES = {Paper No. e67, 43},
      ISSN = {2050-5094},
   MRCLASS = {13D45 (13A35 13A50 13C40 14B15)},
  MRNUMBER = {4623817},
MRREVIEWER = {Adela\ N.\ Vraciu},
       DOI = {10.1017/fms.2023.67},
       URL = {https://doi.org/10.1017/fms.2023.67},
}

@misc{LuoPolstra26,
    author = {Luo, Hongyu and Polstra, Thomas},
     title = {Uniform positivity of {$F$}-signature and effective ideal containments},
      year = {2026},
      note = {Preprint},
}

@incollection {kollar97singsofpairs,
    AUTHOR = {Koll\'ar, J\'anos},
     TITLE = {Singularities of pairs},
 BOOKTITLE = {Algebraic geometry---{S}anta {C}ruz 1995},
    SERIES = {Proc. Sympos. Pure Math.},
    VOLUME = {62, Part 1},
     PAGES = {221--287},
 PUBLISHER = {Amer. Math. Soc., Providence, RI},
      YEAR = {1997},
      ISBN = {0-8218-0894-X; 0-8218-0493-6},
   MRCLASS = {14E30 (14E15 32S40)},
  MRNUMBER = {1492525},
MRREVIEWER = {Alessio\ Corti},
       DOI = {10.1090/pspum/062.1/1492525},
       URL = {https://doi.org/10.1090/pspum/062.1/1492525},
}
\end{document}